\documentclass[11 pt, draft]{article}
\title{Obtaining presentations from group actions without making choices}
\author{Andrew Putman}
\usepackage{amsmath}
\usepackage{amssymb}
\usepackage{amsthm}
\usepackage{epsfig}
\usepackage[margin=1.25in]{geometry}
\usepackage{amsfonts}
\usepackage[font=small,format=plain,labelfont=bf,up,textfont=it,up]{caption}
\usepackage{amscd}
\usepackage{pinlabel}
\usepackage{stmaryrd}
\usepackage{type1cm}
\usepackage{calc}
\usepackage{mathptmx}  

\theoremstyle{plain}
\newtheorem{theorem}{Theorem}[section]

\newtheorem{lemma}[theorem]{Lemma}
\newtheorem{corollary}[theorem]{Corollary}

\newcommand\BeginCases{\setcounter{case}{0}}
\newcommand\BeginSubcases{\setcounter{subcase}{0}}

\theoremstyle{definition}
\newtheorem{definition}[theorem]{Definition}

\newtheorem{case}{Case}
\newtheorem{subcase}{Subcase}

\theoremstyle{remark}
\newtheorem*{remark}{Remark}
\newtheorem*{example}{Example}

\newcommand\Heading[1]{{\bf {\smallskip}{\noindent}#1.}}

\DeclareMathOperator{\Mod}{Mod}
\newcommand\Torelli{\text{${\mathcal I}$}}

\DeclareMathOperator{\SL}{SL}

\newcommand\Bases{\text{$\mathcal{B}$}}

\newcommand\Z{\text{$\mathbb{Z}$}}

\DeclareMathOperator{\HH}{H}

\DeclareMathOperator{\Max}{max}

\newcommand\Th{\text{th}}

\DeclareMathOperator{\Interior}{Interior}

\newcommand\CaptionSpace{\hspace{0.2in}}

\newcommand\Figure[3]{
\begin{figure}[t]
\centering
\centerline{\psfig{file=#2,scale=60}}
\caption{#3}
\label{#1}
\end{figure}}

\newcommand\BigFreeProd{\mathop{\mbox{\Huge{$\ast$}}}}

\newcommand\Move[1]{\text{$\xrightarrow{#1}$}}
\newcommand\Star{\text{Star}}

\begin{document}

\maketitle

\begin{abstract}
Consider a group $G$ acting nicely on a simply-connected
simplicial complex $X$.  Numerous classical methods exist
for using this group action to produce a presentation for $G$.  For
the case that $X/G$ is $2$-connected, we give a new method that
has the novelty that one does not have to identify a fundamental domain
for the action.  Indeed, the resulting presentation is canonical in
the sense that no arbitrary choices need to be made.  It can be viewed
as a nonabelian analogue of a simple result in the study of equivariant homology.
\end{abstract}

\section{Introduction}

A classical theme in group theory is that
if a group $G$ acts nicely on a simply-connected space $X$, then one
can use that action to construct presentations for $G$.  The
investigation of presentations obtained in this way goes back to the $19^{\Th}$ century
study of Fuchsian groups.  One version of such a result (together
with an extensive bibliography) can be found in the paper \cite{BrownPresentation} of Brown,
and this theory has since been subsumed
into the study of Haefliger's theory of complexes of groups \cite{BridsonHaefliger}.

All these classical results require the identification of a fundamental
domain for the action.  While some choice of this type is usually necessary
for the resulting presentation to be finite, there often isn't a
canonical choice of fundamental domain (this is similar to the
fact that most vector spaces do not have a canonical basis).  Moreover,
especially if the action is not cocompact, it can be very difficult to identify
a fundamental domain in a manner explicit enough for the machinery to work.
In this paper, we show how to construct a presentation from a group
action without identifying a fundamental domain, and more generally
without making any arbitrary choices.   

\paragraph{Statement of theorem.}
Let $G$ be a group and $X$ be a simply connected simplicial complex upon which $G$ acts.
We will assume that $G$ acts {\em without rotations}, that is, for
all simplices $s$ of $G$ the stabilizer $G_s$ stabilizes $s$ pointwise (this
can always be arranged to hold by subdividing $X$).
An elementary argument of Armstrong \cite{ArmstrongGenerators} (recalled below
in \S \ref{section:definephi}) shows that if $X/G$ is simply connected, then
$G$ is generated by elements which stabilize vertices.  In other words, we have a surjective
map
$$\BigFreeProd_{v \in X^{(0)}} G_v \longrightarrow G.$$
As notation, if $g \in G$ stabilizes $v \in X^{(0)}$, then we denote $g$
considered as an element of
$$G_v < \BigFreeProd_{v \in X^{(0)}} G_v$$
by $g_v$.  There are then some obvious
elements in the kernel of this map, which we write as relations $f=g$ rather than as elements
$f g^{-1}$.  First, we have $g_v = g_{v'}$ if $v$ and $v'$ are joined
by an edge $e$ and $g \in G_{e}$.  We call these the {\em edge relations}.
Second, we have $g_v h_w g_v^{-1} = (ghg^{-1})_{g(w)}$
for $g \in G_v$ and $h \in G_w$.  We call these relations the {\em conjugation
relations}.  

The following theorem says that if $X/G$ is $2$-connected, then these two families of relations
suffice to give a presentation.
\begin{theorem}[Main theorem]
\label{theorem:presentation}
Let a group $G$ act without rotations on a simply connected simplicial
complex $X$.  Assume that $X/G$ is $2$-connected.  Then
$$G = (\BigFreeProd_{v \in X^{(0)}} G_v)/R,$$
where $R$ is the normal subgroup generated by the conjugation
relations and the edge relations.
\end{theorem}

\begin{remark}
If $X/G$ is $1$-connected but not $2$-connected, then one would also
need relations corresponding to generators for $\pi_2(X/G)$.  One could
extract the precise form of such relations from our proof of Theorem \ref{theorem:presentation}.  We
leave this as an exercise for the interested reader.
\end{remark}

\begin{remark}
The conclusion of Theorem \ref{theorem:presentation} resembles what occurs for groups
acting with {\em strict fundamental domains} (see \cite{BridsonHaefliger}).  Consider a group
$G$ acting without rotations on a $1$-connected simplicial complex $X$.  A subcomplex $C \subset X$
is a strict fundamental domain if it contains exactly one point from each $G$-orbit.  
We then necessarily have $C \cong X/G$.  The theory of complexes of groups shows that
in this situation, the group $G$ can be constructed as
a pushout of the stabilizers of simplices in $C$.  However, Theorem \ref{theorem:presentation}
requires {\em all} of the stabilizer subgroups, and indeed one cannot merely use the stabilizers
of vertices representing the orbits of $G$ (in particular, the assumptions of Theorem \ref{theorem:presentation}
do not imply that $G$ has a strict fundamental domain).  This subtlety already shows up in Armstrong's
theorem mentioned above : one really needs all the vertex stabilizers to generate $G$.  See the remark after
the second example below.
\end{remark}

\paragraph{Examples.}
We now give four examples of group actions to which Theorem \ref{theorem:presentation} can
be applied.  The first two are easy toy examples and the last two are more substantial.

\begin{example}[Amalgamated free products]
If $G = A \ast_C B$,
then Bass-Serre theory \cite{SerreTrees} shows that there is a tree $X$ (i.e.\ a
simply connected $1$-dimensional simplicial complex) upon which $G$ acts
without rotations.  The tree $X$ has the following two properties.
\begin{itemize}
\item $X/G$ is an edge $e$.
\item There exists a lift $\tilde{e}$ of $e$ to $X$ such that if $v,w \in X^{(0)}$ are the
vertices of $\tilde{e}$, then
$$G_{v} = A \quad \text{and} \quad G_w = B \quad \text{and} \quad G_{\tilde{e}} = C.$$
\end{itemize}
The conditions of Theorem \ref{theorem:presentation} are clearly satisfied.  In
the resulting presentation, the stabilizers of vertices correspond to the conjugates
of $A$ and $B$ inside $G$.
\end{example}

\begin{example}[Symmetric group]
Consider the symmetric group $S_n$ on $n$ letters $\{1,\ldots,n\}$.  For distinct $a,b \in \{1,\ldots,n\}$, denote
by $\tau_{a,b}$ the transposition of $a$ and $b$.  We will briefly describe how to use
Theorem \ref{theorem:presentation} to prove that $S_n$ has the following presentation.
\begin{itemize}
\item The generating set $S$ is $\{\text{$\tau_{a,b}$ $|$ $a,b \in \{1,\ldots,n\}$ distinct}\}$.
\item The relations consist of the following for all distinct $a,b \in \{1,\ldots,n\}$ and $s \in S$.
$$\tau_{a,b} = \tau_{b,a} \quad \quad \text{and} \quad \quad \tau_{a,b}^2=1 \quad \quad \text{and} \quad \quad s \tau_{a,b} s^{-1} = \tau_{s(a),s(b)}.$$
\end{itemize}
The proof is by induction on $n$.  The base cases are $n \leq 3$, where the presentation
is easily verified.  Assume now that $n > 3$.

Define $P_n$ to be the following poset.
The elements of $P_n$ are ordered sequences $\sigma = (x_1,\ldots,x_m)$, where the $x_i$
are distinct elements of $\{1,\ldots,n\}$.  For
$\sigma,\sigma' \in P_n$, we have $\sigma \leq \sigma'$ if $\sigma$ is a subsequence
of $\sigma'$.  Next, let $X_n$ be the geometric realization of $P_n$, i.e.\ the simplicial complex whose
$k$-simplices are totally ordered chains
$$\sigma_0 < \sigma_1 < \cdots < \sigma_k$$
of elements of $P_n$.  It is clear that $S_n$ acts on $X_n$ without rotations.  
Bj\"{o}rner and Wachs \cite{BjornerWachs} proved that $X_n$ is $(n-2)$-connected.
It is not hard to show that $X_n/S_n$ is also $(n-2)$-connected (we omit the proof since it is
tangential to the remainder of the paper).  We can thus apply Theorem \ref{theorem:presentation} to the action
of $S_n$ on $X_n$.

If $\sigma = (x_1,\ldots,x_m) \in P_n$, then
$(S_n)_{\sigma}$ is required to fix $\{x_1,\ldots,x_m\}$ pointwise.  It is thus the permutation group
of $\{1,\ldots,n\} \setminus \{x_1,\ldots,x_m\}$, a set with $n-m \leq n-1$ elements.  By induction, these stabilizer
subgroups have presentations of the desired form.  The edge relations identify identical transpositions
that lie in different stabilizer subgroups, and after performing these identifications the conjugation relations 
become the relations $s \tau_{a,b} s^{-1} = \tau_{s(a),s(b)}$ from our presentation.  We conclude that $S_n$ has
a presentation of the desired form.

\begin{remark}
For $k \geq 1$, let $P_n^k$ be the subposet of $P_n$ consisting of ordered sequences $(x_1,\ldots,x_m)$ with
$m \geq k$.  Also, let $X_n^k$ be the geometric realization of $P_n^k$.  The paper of Bj\"{o}rner and Wachs
mentioned in the previous paragraph proves that $X_n^k$ is $(n-k-1)$-connected, and again it is
not hard to show that $X_n^k/S_n$ is $(n-k-1)$-connected as well.
Setting $k = n-2$, the spaces $X_n^k$ and $X_n^k / S_n$ are both $1$-connected, so we can
apply Armstrong's theorem (mentioned before the statement of
Theorem \ref{theorem:presentation}) to get that $S_n$ is generated by stabilizers of
$\sigma \in P_n^k$.  It is clear that stabilizers of length $n-1$ and length $n$ sequences
are trivial.  If $\sigma = (x_1,\ldots,x_k) \in P_n^k$, then the stabilizer
$(S_n)_{\sigma}$ is the permutation group of $\{1,\ldots,n\} \setminus \sigma$, a set with $2$ elements.
We therefore recover (in a maximally complicated manner) the
fact that $S_n$ is generated by transpositions.  Since $S_n$ acts transitively on elements
$(x_1,\ldots,x_k)$ of $P_n^k$ but is not generated by fewer than $n-1$ transpositions, this
illustrates the fact that we really need {\em all} the vertex stabilizers and not just
representatives of each vertex orbit.
\end{remark}
\end{example}

\begin{example}[Torelli group]
In \cite{PutmanInfinite}, the author uses Theorem \ref{theorem:presentation} to obtain
a natural infinite presentation of the Torelli subgroup of the mapping class group of a surface.  Recall
that the mapping class group $\Mod_g$ is the group of homotopy classes of orientation-preserving 
diffeomorphisms of a closed orientable genus $g$ surface $\Sigma_g$.  The Torelli group $\Torelli_g < \Mod_g$
is the kernel of the action of $\Mod_g$ on $\HH_1(\Sigma_g;\Z)$.  The complex used
in \cite{PutmanInfinite} is a variant on the {\em complex of curves}, which is the simplicial
complex whose simplices are sets $\{\gamma_0,\ldots,\gamma_k\}$ of nontrivial homotopy classes
of simple closed curves on $\Sigma_g$ that can be realized disjointly.

\begin{remark}
While it is known from
work of McCool \cite{McCoolMod} (see also \cite{HatcherThurston}) that $\Mod_g$ is finitely
presentable for all $g$, it is not known whether or not $\Torelli_g$ is finitely presentable
for $g \geq 3$.  All that is known is that it is finitely generated for $g \geq 3$ (see \cite{JohnsonFinite})
and is not finitely generated for $g=2$ (see \cite{McCulloughMiller, MessTorelli}).
\end{remark}
\end{example}

\begin{example}[Congruence subgroups]
In \cite{PutmanCongruence}, the author uses Theorem \ref{theorem:presentation} to obtain
natural finite presentations for the level $2$ and $3$ principal congruence subgroups
of $\SL_n(\Z)$.  Recall that the level $L$ principal congruence subgroup $\Gamma_n(L)$ of $\SL_n(\Z)$
is the kernel of the natural map $\SL_n(\Z) \rightarrow \SL_n(\Z/L)$.  Finite presentations
for $\SL_n(\Z)$ were first found by Magnus (see \cite{MilnorKTheory}), but the presentations
in \cite{PutmanCongruence} seem to be the first presentations for $\Gamma_n(L)$ in the literature.
The complex used in \cite{PutmanCongruence} is as follows.  Let $\Bases_n$ be the simplicial
complex whose simplices are sets $\{v_0,\ldots,v_k\}$ of vectors in $\Z^n$ which form part
of a basis for $\Z^n$.  Work of Maazen \cite{MaazenThesis} shows that $\Bases_n$ and
$\Bases_n / \Gamma_n(L)$ are $(n-2)$-connected for $L \geq 2$.  This shows that $\Bases_n$
satisfies the conditions of Theorem \ref{theorem:presentation} for $n \geq 4$.  The
inductive argument starts with the base case $n=2$; for $n=3$, additional $3$-cells are attached to $\Bases_n$
to make the quotient $2$-connected.

\begin{remark}
Of course, Theorem \ref{theorem:presentation} only yields an infinite presentation for
$\Gamma_n(L)$.  Additional techniques are needed to reduce this to a simple finite presentation (this is
where the restriction $L \in \{2,3\}$ is used).
\end{remark}
\end{example}

\paragraph{Comments on proof.}
Though there are numerous methods for obtaining group presentations from group actions, we have
been unable to derive Theorem \ref{theorem:presentation} from any of the ones of which we are
aware.  Rather, our starting point is the theorem of Armstrong discussed above, which
we recall uses the assumption that $X/G$ is simply-connected to show that the
natural map
$$\psi:(\BigFreeProd_{v \in X^{(0)}} G_v)/R \longrightarrow G$$
is surjective.  In fact, Armstrong gives an algorithm (which we call the ``Armstrong construction'')
which takes an element of $G$ and expresses it in terms of vertex stabilizers.  Numerous
choices need to be made when running this algorithm.  However, via a careful analysis
of the combinatorics of homotoping loops and discs inside simplicial complexes we will show that modulo our
relations none of these choices matter.  The upshot is that the Armstrong construction
induces a well defined surjective map
$$\psi^{-1} : G \longrightarrow (\BigFreeProd_{v \in X^{(0)}} G_v)/R$$
satisfying $\psi \circ \psi^{-1}=1$.

\begin{remark}
Our proof is inspired in part by another paper of Armstrong \cite{ArmstrongBrown} in which
he gives a combinatorial-topological derivation of the presentation of Brown \cite{BrownPresentation}
alluded to in the first paragraph of this paper.  Brown originally derived his
presentation via Bass-Serre theory.
\end{remark} 

\paragraph{Motivation from equivariant topology.}
Theorem \ref{theorem:presentation} can be viewed as a nonabelian analogue
of a simple result in equivariant topology.  We proceed to sketch this.  Nothing
in this section is needed in the rest of the paper.

Let $G$ and $X$ be as in Theorem \ref{theorem:presentation}.  We will assume that $X/G$ is
the simplicial complex whose simplices are orbits of simplices of $X$ (this can always be
arranged by subdividing $X$; we remark that it does not follow from the fact that $G$ acts
without rotations).  Let $Y$ be the $1$-skeleton of $X/G$.  Regarding the simplicial complex $Y$ as a category
whose objects are simplices and whose morphisms are inclusions, there is a natural contravariant functor
$\mathcal{H}_1$ from $Y$ to the category of abelian groups.  Namely, if $s$ is a simplex of $Y$, then
$\mathcal{H}_1(s) = \HH_1(G_{\tilde{s}};\Z)$, where $\tilde{s}$ is a lift of $s$ to $X$.  Define
$\mathcal{C}(\mathcal{H}_1)$ to be the colimit of $\mathcal{H}_1$; i.e.\ the quotient of the abelian group
$$\bigoplus_{\text{$s$ a simplex of $Y$}} \mathcal{H}_1(s)$$
which for all edges $e = \{v,v'\}$ of $Y$ identifies $\mathcal{H}_1(e)$ with its images in 
$\mathcal{H}_1(v)$ and $\mathcal{H}_1(v')$.

Let $\HH_k^{G}(X;\Z)$ denote the equivariant homology groups of $G$ acting on $X$ in the sense
of \cite[\S VII.7]{BrownCohomology} (defined in terms of the Borel construction).  These satisfy the following two properties.
\begin{itemize}
\item Since $\pi_1(X)=1$, the spectral sequence (7.2) of \cite[\S VII.7]{BrownCohomology}
implies that $\HH_1^{G}(X;\Z) \cong \HH_1(G;\Z)$.
\item The spectral sequence whose $E^2$ page is described in \cite[\S VII.8]{BrownCohomology} is a 
first-quadrant spectral sequence, so it induces a 5-term exact sequence.  This exact sequence takes
the form
\begin{equation}
\label{eqn:5term}
\HH_2^{G}(X;\Z) \longrightarrow \HH_2(X/G;\Z) \longrightarrow \mathcal{C}(\mathcal{H}_1) \longrightarrow \HH_1^{G}(X;\Z) \longrightarrow \HH_1(X/G;\Z) \longrightarrow 0.
\end{equation}
\end{itemize}
Since $X/G$ is 2-connected, \eqref{eqn:5term} implies that
\begin{equation}
\label{eqn:equivarident}
\HH_1(G;\Z) \cong \HH_1^{G}(X;\Z) \cong \mathcal{C}(\mathcal{H}_1).
\end{equation}
One can view the group presentation in Theorem \ref{theorem:presentation} as a nonabelian analogue
of $\mathcal{C}(\mathcal{H}_1)$, 
and in fact \eqref{eqn:equivarident} can be easily deduced from Theorem \ref{theorem:presentation}.

\begin{remark}
In fact, our original motivation for guessing that something like Theorem \ref{theorem:presentation}
might be true was van den Berg's unpublished thesis \cite{VanDenBerg}, where she used
\eqref{eqn:equivarident} to give a new proof of Johnson's theorem \cite{JohnsonAbel} giving the abelianization of the
Torelli subgroup of the mapping class group.
\end{remark}

\paragraph{Acknowledgments.}
This paper is a revised version of part of my thesis, and I wish to thank my advisor Benson
Farb for his help and encouragement.  Additionally, I wish to thank Anne Thomas, Shmuel Weinberger, 
and Ben Wieland for helpful conversations.

\section{Simplicial Preliminaries}

The heart of our proof is a careful examination of the combinatorics
of homotoping loops and discs in simplicial complexes.  In this section,
we establish some preliminary results in this direction.

\subsection{Simplicial complexes}

We first establish our notation for simplicial complexes.  Let $X$ be a simplicial
complex.  We will denote by $|X|$ the geometric realization of $X$.  When we say
that a set $\sigma = \{x_0,\ldots,x_k\}$ of vertices of $X$ forms a $k$-simplex, we allow the
possibility that $x_i = x_j$ for some $i$ and $j$.  By a path or a loop
in $X$, we mean a simplicial path or loop in the $1$-skeleton.  We will denote the
path or loop that starts at a vertex $v_0$, goes along an edge to a vertex $v_1$, then goes
along an edge to a vertex $v_2$, etc.\ and
ends at a vertex $v_n$ (which equals $v_0$ if the path is a loop) by
$v_0 - v_1 - \cdots - v_n$.  To simplify our notation, we will allow the possibility
that $v_i = v_{i+1}$ for some $0 \leq i < n$ (this is consistent with regarding
$\{v,v\}$ as a degenerate edge for a vertex $v$ of $X$).  However, we will not regard
such a path as injective (for example, in the definition of degenerate discs below 
in \S \ref{section:loopmoves}).
If $f : X \rightarrow Y$ is a map of simplicial complexes and $\gamma$ is a path or loop in $X$, then
we will denote by $f_{\ast}(\gamma)$ the induced path or loop in $Y$.  Finally, if $s$ is a simplex
of $X$, then the {\em star} of $s$, denoted $\Star_X(s)$, is the subcomplex of $X$ consisting of all
simplices $s'$ such that there is a simplex $s''$ of $X$ containing both $s$ and $s'$ as faces.

\subsection{Collapsing loops via degenerate discs}
\label{section:loopmoves}

While simplicially contracting paths to points, we will need the notion of a degenerate disc.  To define
this, we need a pair of preliminary definitions.

\begin{definition}
A {\em space/loop pair} is a pair $(X,\gamma)$, where $X$ is a simplicial complex with a basepoint $\ast$ 
and $\gamma$ is a loop in $X$ that is based at $\ast$.
\end{definition}

\begin{definition}
A {\em nondegenerate disc} is a space/loop pair $(D,\gamma)$ with the following properties.
\begin{itemize}
\item $D$ is homeomorphic to a closed 2-disc.
\item The basepoint $\ast \in D$ lies in $\partial D$.
\item $\gamma$ is an injective loop in $\partial D$ that is based at $\ast$ and goes once
around $\partial D$.
\end{itemize}
\end{definition}

\begin{definition}
A {\em degenerate disc} (see Figure \ref{figure:degeneratedisc}.b) is a space/loop pair $(D,\gamma)$ such that
there exists a sequence $(D_1,\gamma_1),\ldots,(D_k,\gamma_k)$ of space/loop pairs with the following properties.
\begin{itemize}
\item $(D_1,\gamma_1)$ is a nondegenerate disc and $(D_k,\gamma_k) = (D,\gamma)$.
\item For $1 \leq i < k$, the space/loop pair $(D_{i+1},\gamma_{i+1})$ is obtained from $(D_i,\gamma_i)$ by one
of the following two moves.  Let $\ast$ be the basepoint of $D_i$.
\begin{itemize}
\item If $\gamma_i$ contains a subpath of the form $x_1 - x_2$ for vertices
$x_1,x_2 \in D_i^{(0)}$ and there is some vertex $y \in D_i^{(0)}$ such that
$\{x_1,x_2,y\}$ is a $2$-simplex of $D_i$ and $y \neq x_1,x_2$ (see Figure
\ref{figure:degeneratedisc}.a), then we can homotope $x_1 - x_2$ to
$x_1 - y - x_2$ and delete $\{x_1,x_2\}$ and $\{x_1,x_2,y\}$ from $D_i$.  We will call this a {\em two-dimensional collapse}.
\item If $\gamma_i$ contains a subpath of the form $x_1 - x_2 - x_1$ for vertices
$x_1,x_2 \in D_i^{(0)}$ with $x_2 \neq x_1,\ast$ (see
Figure \ref{figure:degeneratedisc}.a),
then we can homotope $x_1 - x_2 - x_1$ to the constant
path $x_1$ and delete $\{x_1,x_2\}$ from $D_i$.  We will call this a {\em one-dimensional collapse}.
\end{itemize}
We will use the term {\em collapse} to refer to either a one- or two-dimensional collapse.
\end{itemize}
\end{definition}

\Figure{figure:degeneratedisc}{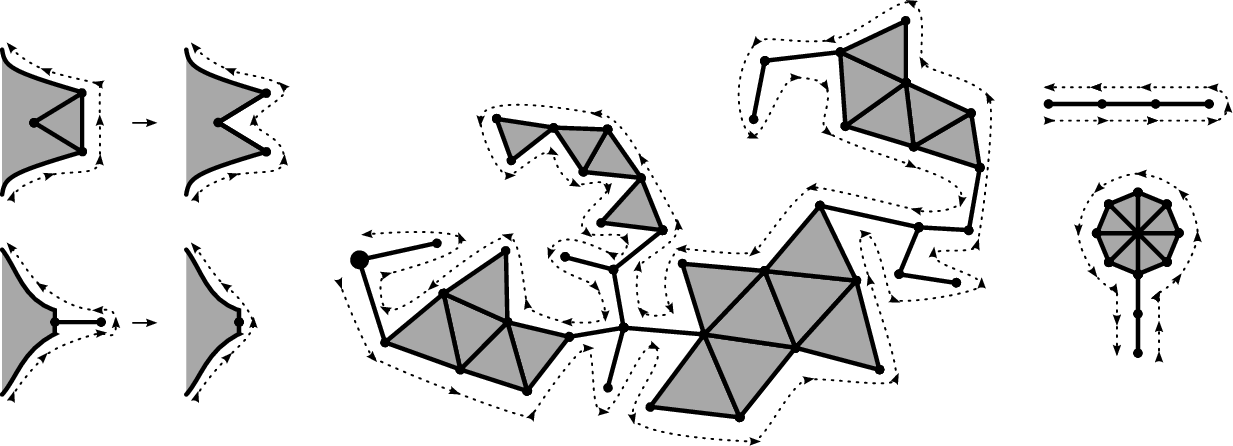}{a. Two- and one-dimensional collapses \CaptionSpace 
\CaptionSpace b. A degenerate disc \CaptionSpace c. A line and a flower}

\begin{remark}
It is easy to see that the one- and two-dimensional collapses involved in the definition of a degenerate disc do
not delete simplices that appear elsewhere in the loop in question.
\end{remark}

\begin{remark}
The intermediate space/loop pairs $(D_i,\gamma_i)$ in the definition of a degenerate discs are themselves degenerate discs, but
need not be nondegenerate discs.
\end{remark}

\begin{remark}
The cell complexes used in van Kampen diagrams are (almost) degenerate discs, the only difference
being that we require each cell to be a triangle.
\end{remark}

A key technical lemma concerning degenerate discs is as follows.

\begin{lemma}
\label{lemma:trianglelemma}
Let $(D,\gamma)$ be a degenerate disc and let $X$ be a subcomplex of $D$ which is homeomorphic
to a closed 2-disc.  Let $e$ be an edge of $\partial X$ that does not lie in $\gamma$.  There then exists a sequence
$s_1,\ldots,s_k$ of $2$-simplices of $D$ with the following properties.
\begin{itemize}
\item $s_i$ is not a $2$-simplex of $X$ for $1 \leq i \leq k$.
\item $s_i \cap s_{i+1}$ is an edge of $D$ for $1 \leq i < k$.
\item $e$ is an edge of $s_1$ and $s_k$ contains an edge of $\gamma$.
\end{itemize}
\end{lemma}
\begin{proof}
Assume first that $(D,\gamma)$ is a nondegenerate disc.  Observe that each component of $|D| \setminus |X|$ must
contain an edge of $\gamma$.  We can therefore choose a continuous path $\delta$ in $|D| \setminus \Interior(|X|)$ 
which begins at an interior point of $e$ and ends at an interior point of an edge of $\gamma$.  Moreover, we can choose $\delta$
such that it does not pass through any vertices of $D$, such that all of its intersections with edges of $D$ are transverse, 
and such that it only intersects finitely many $2$-simplices of $D$.  The desired sequence of $2$-simplices 
is then the sequence of $2$-simplices through which $\delta$ passes.

We now consider the general case.  It is enough to show that if the lemma is true for some degenerate disc
$(D',\gamma')$, then it remains true for the result $(D,\gamma)$ of performing a
collapse to $(D',\gamma')$.  The case of a one-dimensional collapse being trivial, we consider the
case that $(D,\gamma)$ is the result of performing a two-dimensional collapse to $(D',\gamma')$.  Let $t$
be the $2$-simplex of $D'$ that is collapsed.  Consider
a subcomplex $X$ and an edge $e$ of $D$ as in the lemma.  Regarding $D$ and thus $X$ as a subcomplex of $D'$,
by assumption there is a sequence $s_1,\ldots,s_k$ of $2$-simplices of $D'$ satisfying
the conditions of the lemma (applied to the subcomplex $X$ in the degenerate disc $(D',\gamma')$).  
Observe that since $e$ is not an edge of $\gamma$, it follows that $s_1 \neq t$.  Set 
$$l = \Max\{\text{$i$ $|$ $1 \leq i \leq k$, $s_j \neq t$ for $1 \leq j \leq i$}\}.$$
It is then clear that the sequence $s_1,\ldots,s_l$ of $2$-simplices of $D$ satisfies the conditions of the lemma,
and we are done.
\end{proof}

Using Lemma \ref{lemma:trianglelemma}, we can show that degenerate discs can always be contracted
to points by a sequence of collapses.  This is a standard result, but we give the short proof
to make this paper self-contained.  The reader should be warned that analogous statements are false in dimensions
greater than $2$ (see \cite[\S 4, Example 3]{BingSurvey} and \cite{Goodrick}).

\begin{corollary}
\label{corollary:cancollapse}
Let $(D,\gamma)$ be a degenerate disc.  Then we can perform a sequence of collapses to $(D,\gamma)$
such that the resulting degenerate disc $(D',\gamma')$ consists of a single vertex $\ast$ together with
the constant path $\gamma' = \ast$.
\end{corollary}
\begin{proof}
If $D$ contains any $2$-simplex, then by Lemma \ref{lemma:trianglelemma} there is some $2$-simplex containing
an edge of $\gamma$.  We can therefore perform a sequence of two-dimensional collapses so as to obtain
a degenerate disc $(D',\gamma')$ such that $D'$ is a simply-connected $1$-complex; i.e.\ a tree.  We
can then perform a sequence of one-dimensional collapses to collapse $D'$ to $\ast$, and we are done.
\end{proof}

Finally, the following two special types of degenerate discs will play a key role in our proofs.

\begin{definition}
A {\em line} (see top of Figure \ref{figure:degeneratedisc}.d) is a degenerate disc $(D,\gamma)$ of
the following form.
\begin{itemize}
\item $D$ is the complex consisting of vertices $v_1,\ldots,v_n$ for some $n \geq 1$ together with edges between
$v_i$ and $v_{i+1}$ for $1 \leq i < n$.  The basepoint $\ast$ is $v_1$.  If $n \geq 2$, then
the edge $v_{n-1} - v_n$ will be called the {\em final edge}.
\item $\gamma$ is the path $v_1 - \cdots - v_{n-1} - v_n - v_{n-1} - \cdots - v_1$.
\end{itemize}
A {\em flower} (see bottom of Figure \ref{figure:degeneratedisc}.d) is a degenerate disc $(D,\gamma)$ of the
following form.
\begin{itemize}
\item There is a nondegenerate disc $(D',\gamma')$ satisfying the following condition.  
The complex $D'$ is a subcomplex of
$D$ and for some $n \geq 1$ the complex $D'$ consists of $D$, vertices $v_1,\ldots,v_{n-1}$, and
edges between $v_i$ and $v_{i+1}$ for $1 \leq i < n$, where $v_n$ is the basepoint of $(D',\gamma')$.  The
basepoint of $(D,\gamma)$ is $v_1$.
\item $\gamma'$ is the path $v_1 - \cdots - v_n - \gamma - v_n - \cdots - v_1$.
\end{itemize}
We will call $(D',\gamma')$ the {\em bloom} of $(D,\gamma)$ and the path $v_1 - \cdots - v_n$ the {\em stem}; we
will confuse the stem with its associated subcomplex.
\end{definition}

\begin{remark}
A line may consist of a single vertex.  Similarly, the stem of a flower may consist of a 
single vertex, in which case the bloom is the entire flower.
\end{remark}

The following corollary to Lemma \ref{lemma:trianglelemma} will be frequently used, and is the key
reason we introduced lines and flowers.

\begin{corollary}
\label{corollary:lineflower}
Let $(D,\gamma)$ be a degenerate disc.
\begin{enumerate}
\item Let $e$ be an edge of $D$.  Let $(D',\gamma')$ be a minimal degenerate disc with $e$ an edge of $D'$ 
to which we can collapse $(D,\gamma)$.  Then $(D',\gamma')$ is a line whose final edge is $e$.
\item Let $X$ be a subcomplex of $D$ which is homeomorphic to a closed 2-disc.  Let $(D',\gamma')$ be a minimal degenerate
disc with $X$ a subcomplex of $D'$ to which we can collapse $(D,\gamma)$.  Then
$(D',\gamma')$ is a flower whose bloom consists of $(X,\gamma'')$ for some path $\gamma''$
around the boundary of $X$.  If $X$ contains the basepoint $\ast$ of $(D,\gamma)$, then in fact $D' = X$.
\end{enumerate}
\end{corollary}
\begin{proof}
The proofs of the two conclusions are similar; we will prove the more difficult second one and leave the first to the reader.  

The first step is to prove that any $2$-simplex of $D'$ is contained in $X$.  Assume otherwise.  If $\partial X$ contains
an edge $e$ which is not contained in $\gamma'$, then by Lemma \ref{lemma:trianglelemma} (applied to
$X$ and $e$), there must exist some $2$-simplex $s$ such that $s$ is not a $2$-simplex of $X$ and such that 
$s$ contains an edge of $\gamma'$.  We can then perform a two-dimensional collapse so as to eliminate $s$, 
contradicting the minimality of $(D',\gamma')$.  All the edges of $\partial X$ must therefore be contained
in $\gamma'$.  

An easy induction establishes that if $(E,\delta)$ is any degenerate disc and $\varepsilon$ is an
edge of $\delta$, then at most one $2$-simplex of $E$ contains $\varepsilon$.  We conclude that there 
does not exist a $2$-simplex of $D'$ whose intersection with $X$
is an edge.  Applying Lemma \ref{lemma:trianglelemma} again, we conclude that since there exists a 
$2$-simplex $t$ of $D'$ that does not lie in $X$, there must exist a $2$-simplex $t'$ of $D'$ that does not lie in 
$X$ such that $t'$ 
contains an edge of $\gamma'$.
We can therefore perform a two-dimensional collapse so as to eliminate $t'$, 
again contradicting the minimality of $(D',\gamma')$.  We conclude that every $2$-simplex of $D'$
is in fact contained in $X$.

Since $D'$ and $X$ are contractible, the space $D'/X$ is contractible.  Since $D'/X$ contains no $2$-cells, 
it follows that $D'/X$ is a tree.  We conclude that the closure of every component of $|D'| \setminus |X|$ must be a 
tree that intersects $X$ in exactly one point.  Since trees
can be collapsed to any of their vertices, we deduce that $|D'| \setminus |X|$ can have at most one component (the
one containing the basepoint of $(D',\gamma')$; this exists exactly when $X$ does not contain the basepoint
of $(D',\gamma')$).  Moreover, if there is such a component, then its closure must satisfy the conditions
of the stem of a flower, and we are done.
\end{proof}

\subsection{Simplicial homotopy}
\label{section:simplicialhomotopy}

We now give moves by which we will homotope discs in simplicial complexes.
The following definition is a variant on the notion of {\em contiguity classes} in \cite{Spanier}.

\begin{definition}
Let $(K,L)$ be a simplicial pair, let $X$ be a simplicial complex, and let $f:K \rightarrow X$ be a simplicial
map.
\begin{itemize}
\item Let $s$ be an $n$-simplex of $K$ with $n \geq 1$ and let $v$ be a vertex of $s$.  Let $(K',L')$ be the
result of subdividing the simplicial pair $(K,L)$ by adding a vertex $v'$ to the interior of $s$, and define
$f':K' \rightarrow X$ by setting
$$f'(w) = 
\begin{cases}
f(w) & \text{if $w \neq v'$}\\
f(v) & \text{if $w = v'$}
\end{cases} \quad \quad \quad (w \in (K')^{(0)})$$
and extending linearly.  We will say that the complex $(K',L')$ and the map $f':K' \rightarrow X$ come from performing
an {\em elementary subdivision} of $f$ along $s$ with $v$.
\item Let $v \in K^{(0)} \setminus L^{(0)}$ and $x \in X^{(0)}$ 
be such that for all simplices $s$ of $\Star_K(v)$,
the set $f(s) \cup \{x\}$ is a simplex of $X$.  Define $f':K \rightarrow X$ by setting
$$f'(w) =
\begin{cases}
f(w) & \text{if $w \neq v$}\\
x & \text{if $w = v$}
\end{cases} \quad \quad \quad (w \in K^{(0)})$$
and extending linearly.  We will say that the map $f':K \rightarrow X$ comes from performing an {\em elementary
push} of $f$ along $v$ with $x$.
\item If $f':K \rightarrow X$ differs from $f$ by a sequence of elementary pushes, then we
say that $f$ and $f'$ are {\em elementarily equivalent}.
\end{itemize}
\end{definition}

We will need the following theorem.  It is essentially \cite[Theorem 3.5.6]{Spanier}, but since our definitions
are a little different we include a sketch of the proof.

\begin{theorem}[{\cite[Theorem 3.5.6]{Spanier}}]
\label{theorem:elementaryequiv}
Let $(K,L)$ be a compact simplicial pair, let $X$ be a simplicial complex, and let $f_0,f_1:K \rightarrow X$
be simplicial maps satisfying the following conditions.
\begin{itemize}
\item $f_0|_{L} = f_1|_{L}$.
\item The maps $f_0$ and $f_1$ are homotopic fixing $L$; i.e.\ there is a continuous map
$F : |K| \times [0,1] \rightarrow |X|$ such that $F(\cdot,0) = f_0$, such that $F(\cdot,1) = f_1$, and such
that $F(x,t) = f_0(x)$ for $x \in |L|$.
\end{itemize}
Then there exists a subdivision $(K',L')$ of $(K,L)$ and simplicial maps
$f_1',f_2':K \rightarrow X$ which are obtained from the $f_i$ by sequences of elementary subdivisions such
that $f_1'$ and $f_2'$ are elementarily equivalent.
\end{theorem}
\begin{proof}
Since $K$ is compact, a Lebesgue number argument implies that there exists a sequence
$0 = t_0 < t_1 < \cdots < t_n = 1$ such that for $x \in |K|$ and $1 \leq i \leq n$, there
exists a vertex $v \in X^{(0)}$ such that $F(x,t_{i-1}), F(x,t_i) \in \Interior(\Star_X(v))$.  For
$0 \leq i \leq n$, define $g_i : |K| \rightarrow |X|$ by $g_i(x) = F(x,t_i)$ for $x \in |K|$, so
$g_0 = f_0$ and $g_n = f_1$.  

By construction, for $1 \leq i \leq n$ the set
$$\mathcal{U}_i := \{\text{$g_{i-1}^{-1}(\Interior(\Star_X(v))) \cap g_i^{-1}(\Interior(\Star_X(v)))$ $|$ $v \in X^{(0)}$}\}$$
is an open cover of $K$.  Let $(K',L')$ be a subdivision of $(K,L)$ which is finer than $\mathcal{U}_i$ for
$1 \leq i \leq n$.  For $0 \leq i \leq n$, we can construct simplicial maps $\phi_i : K' \rightarrow X$ satisfying
the following conditions.
\begin{enumerate}
\item $\phi_0$ and $\phi_n$ can be obtained by a sequence of elementary subdivisions from $f_0$ and $f_1$, respectively.
\item For all vertices $v \in (K')^{(0)}$, we have $f_{i-1}(\Star_{K'}(v)) \cup f_i(\Star_{K'}(v)) \subset \Star_X(\phi_i(v))$
for $1 \leq i \leq n$.
\item For all vertices $v \in (L')^{(0)}$, we have $\phi_i(v) = \phi_j(v)$ for all $0 \leq i,j \leq n$.
\end{enumerate}
For $1 \leq i \leq n$, items 2 and 3 imply that $\phi_i$ is elementarily equivalent to $\phi_{i-1}$.  We conclude
that $\phi_0$ is elementarily equivalent to $\phi_n$, so $f_1' = \phi_0$ and $f_2' = \phi_n$ satisfy
the conditions of the theorem.
\end{proof}

\section{Proof of the main theorem}
\label{section:proof}

We begin by observing that since $G$ acts on $X$ without rotations, we can subdivide $X$
without affecting the conclusion of the theorem.  Indeed, subdividing introduces new vertices, but the stabilizer $G_{x'}$
of a new vertex $x'$ lies in $G_x$ for some old vertex $x$, and the edge relations identify $G_{x'}$ with its
image in $G_x$.  Moreover, it is easy to see that all the new edge and conjugation relations involving elements of $G_{x'}$
are consequences of the old edge and conjugation relations.
By taking the barycentric subdivision, we can assure that $X/G$ is the simplicial complex whose simplices are
the orbits of simplices in $X$.

\begin{remark}
If $G$ did not act without rotations, then we would need to take the second barycentric subdivision to assure
that $X/G$ is the simplicial complex whose simplices are the orbits of simplices in $X$.
\end{remark}

Let $\pi : X \rightarrow X/G$ be the projection and let
$$\Gamma = (\BigFreeProd_{x \in X^{(0)}} G_x)/R$$
be as in the statement of the theorem.  As in the introduction, for $x \in X^{(0)}$ and $h \in G_x$, we will
denote by $h_x \in \Gamma$ the corresponding element of $G_x < \Gamma$.  There is an obvious homomorphism 
$\psi : \Gamma \rightarrow G$.  We
will construct a surjective homomorphism $\psi^{-1}: G \rightarrow \Gamma$ such
that $\psi \circ \psi^{-1}=1$; the theorem will immediately follow.

The proof will have three parts.  In \S \ref{section:definephi}, we give a procedure (due to
Armstrong \cite{ArmstrongGenerators}) for expressing an element of $G$
as a product of vertex stabilizers.  This procedure is reminiscent of standard
arguments involving covering spaces.  In \S \ref{section:phiwelldefined}, we will show that the
resulting element of $\Gamma$ is independent of the choices made.  This will
define $\psi^{-1}$.  Finally, in \S \ref{section:phihomo} we will show that $\psi^{-1}$ is a surjective
homomorphism.

Throughout the whole proof, we will fix some basepoint $\tilde{v} \in X^{(0)}$ and define $v = \pi(\tilde{v}) \in X/G$.

\subsection{Expressing elements of $\boldsymbol{G}$ as products of vertex stabilizers : the Armstrong construction}
\label{section:definephi}

Consider $g \in G$.  Let $\tilde{\gamma}$ be a simplicial path in $X$ from $\tilde{v}$ to $g(\tilde{v})$.
The projection $\pi$ sends $\tilde{\gamma}$ to a closed loop in $X/G$ based at $v$.
Since $X/G$ is simply-connected, there is some simplicial map $\phi : (D,\ast) \rightarrow (X/G,v)$, where
$(D,\gamma)$ is a nondegenerate disc and $\phi_{\ast}(\gamma) = \pi_{\ast}(\tilde{\gamma})$ (in fact, 
for later use we point out that in our construction, we will never use the nondegeneracy of $(D,\gamma)$).  This
is all illustrated in Figure \ref{figure:defineinverse}.a. 

Since $(D,\gamma)$ is a degenerate disc, Corollary \ref{corollary:cancollapse} says that there is a sequence
$$(D,\gamma) = (D_1,\gamma_1), (D_2,\gamma_2), \ldots, (D_n,\gamma_n) = (\ast,\ast)$$
of degenerate discs such that for $1 \leq i < n$, the degenerate disc
$(D_{i+1},\gamma_{i+1})$ differs from $(D_i,\gamma_i)$ by either a one- or two-dimensional collapse.  
Each $D_i$ is a subcomplex of $D$, so $\phi$ restricts to a map $\phi_i : D_i \rightarrow X/G$.  The resulting 
sequence of based loops $(\phi_i)_{\ast}(\gamma_i)$ gives a simplicial homotopy from $\pi_{\ast}(\tilde{\gamma})$
to the constant path.  We will construct a sequence of paths $\tilde{\gamma}_1,\ldots,\tilde{\gamma}_n$ in $X$
which all start at $v$ and which satisfy $\pi_{\ast}(\tilde{\gamma}_i) = (\phi_i)_{\ast}(\gamma_i)$ for $1 \leq i \leq n$.
Additionally, our construction will yield a sequence $\tilde{v}_1,\ldots,\tilde{v}_n$ of vertices of $X$
and a sequence $h_1,\ldots,h_n$ of elements of $G$ such that $h_i \in G_{\tilde{v}_i}$ for $1 \leq i \leq n$
and such that $h_1^{-1} \cdots h_n^{-1} = g$.

\Figure{figure:defineinverse}{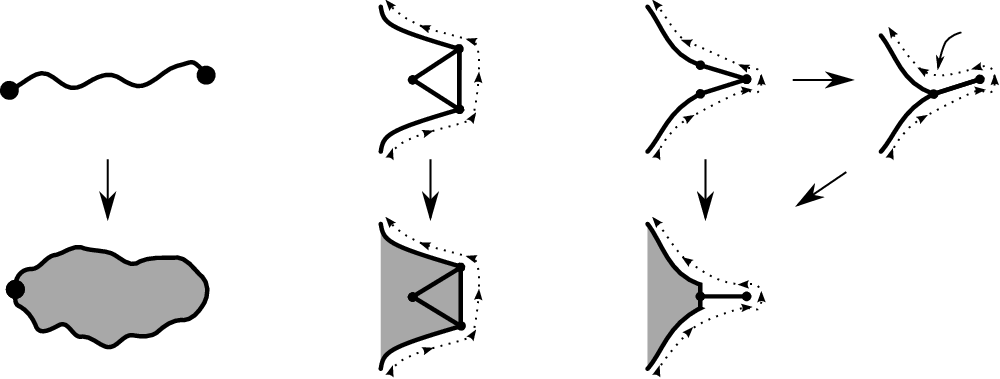}{a. $\tilde{\gamma}$ and its projection \CaptionSpace b. Lifting
a two-dimensional collapse \CaptionSpace c. Swinging around a pivot point to lift a one-dimensional collapse}

Begin by defining $\tilde{\gamma}_1 = \tilde{\gamma}$.  Assume
that for some $1 \leq k < n$ we have constructed paths $\tilde{\gamma}_1,\ldots,\tilde{\gamma}_k$ with
the indicated properties.  If $(D_{k+1},\gamma_{k+1})$ differs from $(D_{k},\gamma_{k})$
by a two-dimensional collapse, then by definition
there is some simplex $\{x_1,x_2,y\}$ of $D_k$ such that $x_1 - x_2$ is a
subpath of $\gamma_k$ which is homotoped to the subpath $x_1 - y - x_2$ of
$\gamma_{k+1}$ (see the bottom portion of Figure \ref{figure:defineinverse}.b).
Let $\tilde{x}_1 - \tilde{x}_2$ be the relevant portion of $\tilde{\gamma}$.  
Since simplices in $X/G$ are exactly the orbits of simplices in $X$, there
exists some $\tilde{y} \in X^{(0)}$
such that $\pi(\tilde{y})=\phi_k(y)$ and such that $\{\tilde{x}_1,\tilde{x}_2, \tilde{y}\}$
is a simplex of $X$ (see the top portion of Figure \ref{figure:defineinverse}.b; we wish
to point out that since $\phi_k$ need not be injective, we might have $\tilde{y} = \tilde{x}_i$
for some $i$).
Let $\tilde{\gamma}_{k+1}$ equal $\tilde{\gamma}_k$ with the subpath
$\tilde{x}_1 - \tilde{x}_2$ replaced with $\tilde{x}_1 - \tilde{y} - \tilde{x}_2$.  It is
clear that $\tilde{\gamma}_{k+1}$ projects to $(\phi_{k+1})_{\ast}(\gamma_{k+1})$.  In this
case, we set $\tilde{v}_k = \tilde{v}$ and $h_k = 1$.

\begin{remark}
Since $h_k=1$, the precise value of $\tilde{v}_k$ is immaterial.
\end{remark}

If $(D_{k+1},\gamma_{k+1})$ differs from $(D_{k},\gamma_{k})$ by a one-dimensional collapse,
then we may have to modify $\tilde{\gamma}_k$ before homotoping it.
Indeed, let $x_1 - x_2 - x_1$ be the relevant subpath of $\gamma_{k}$.
As indicated in Figure \ref{figure:defineinverse}.c, the portion of $\tilde{\gamma}_k$
which projects to $\phi_k(x_1) - \phi_k(x_2) - \phi_k(x_1)$ may be of the form $\tilde{x}_1 - \tilde{x}_2 - \tilde{x}_1'$
with $\pi(\tilde{x}_1) = \pi(\tilde{x}_1')$ but with $\tilde{x}_1 \neq \tilde{x}_1'$ (observe that
if this happens, then we must have $\phi_k(x_1) \neq \phi_k(x_2)$).  However, since simplices
of $X/G$ are the orbits of simplices in $X$, there must exist some $h \in G_{\tilde{x}_2}$ such that 
$h_k(\tilde{x}_1') = \tilde{x}_1$ (if $\tilde{x}_1' = \tilde{x}_1$, then $h = 1$).  Let $\rho_1$
be the portion of $\tilde{\gamma}_k$ before $\tilde{x}_1 - \tilde{x}_2 - \tilde{x}_1'$ and $\rho_2$ be
the portion after it.  We then define $\tilde{\gamma}_{k+1}'$ to equal
$$\rho_1 - \tilde{x}_1 - \tilde{x}_2 - \tilde{x}_1 - h_k(\rho_2);$$
see the top right hand portion of Figure \ref{figure:defineinverse}.c.  
We will call this the result of {\em swinging $\tilde{\gamma}_k$ around the pivot point $\tilde{x}_2$ by $h_k$}.  We
then can remove $\tilde{x}_2$ from $\tilde{\gamma}_{k+1}'$ to get $\tilde{\gamma}_{k+1}$, a lift
of $(\phi_{k+1})_{\ast}(\gamma_{k+1})$.  Define $\tilde{v}_k = \tilde{x}_2$ and $h_k = h$.

\begin{remark}
Technically speaking, our pivot point $\tilde{x}_2$ lies in $X$, not on $\tilde{\gamma}$, and in fact $\tilde{\gamma}$
may pass through $\tilde{x}_2$ multiple times.  However, to simplify our notation we will abuse notation and always assume that
pivot points correspond to specific points on our paths.  In particular, if we say that two pivot points on a path
are the same, we mean not merely that they are the same point of $X$ but that moreover they correspond to the
same point on the path.
\end{remark}

We have now constructed the paths $\tilde{\gamma}_1,\ldots,\tilde{\gamma}_n$, the vertices
$\tilde{v}_1,\ldots,\tilde{v}_{n-1}$, and the group elements $h_1,\ldots,h_{n-1}$ with $h_i \in G_{\tilde{v}_i}$
for $1 \leq i \leq n-1$.  For $1 \leq i < n$, it is immediate from our construction that
the endpoint of $\tilde{\gamma}_{i+1}$ is equal to the result of applying $h_i$ to the endpoint of
$\tilde{\gamma}_i$.  Since the endpoint of $\tilde{\gamma}_1$ is $g(\tilde{v})$, we deduce that the
endpoint of $\tilde{\gamma}_n$ is $h_{n-1} \cdots h_1 \cdot g(\tilde{v})$.  But $\tilde{\gamma}_n$
is the constant path $\tilde{v}$.  We conclude that 
$$h_{n-1} \cdots h_1 \cdot g(\tilde{v}) = \tilde{v}.$$
Define $\tilde{v}_n = \tilde{v}$ and $h_n = (h_{n-1} \cdots h_1 \cdot g)^{-1}$.  We thus have
$h_n \in G_{\tilde{v}_n}$ and
$$g = h_1^{-1} \cdots h_n^{-1},$$
as desired.  

Observe that  
$$\psi((h_1)^{-1}_{\tilde{v}_1} \cdots (h_n)^{-1}_{\tilde{v}_n}) = g.$$
We will say that
$$(h_1)^{-1}_{\tilde{v}_1} \cdots (h_n)^{-1}_{\tilde{v}_n} \in \Gamma$$
is the result of applying the {\em Armstrong construction} to $g$.  Though 
a priori the output of the Armstrong construction depends on numerous choices, in \S \ref{section:phiwelldefined}
we will show that it is in fact independent of those choices.

For later use, we now define some notation.  
A degenerate disc $(D',\gamma')$ together with a lift $\tilde{\gamma}'$ to $X$
of the image of $\gamma'$ under some (implied) map from $D'$ to $X/G$ will be denoted by $(D',\gamma',\tilde{\gamma}')$.
We will denote the transition from $(D',\gamma',\tilde{\gamma}')$ to $(D'',\gamma'',\tilde{\gamma}'')$ coming
from a one- or two-dimensional collapse plus swinging $\tilde{\gamma}'$ around the pivot point $\tilde{v}'$ by 
$h' \in G_{\tilde{x}'}$ by
$$(D',\gamma',\tilde{\gamma}') \Move{h',\tilde{v}'} (D'',\gamma'',\tilde{\gamma}'').$$
We will denote the sequence of moves given by the above construction
by
\begin{equation}
\label{eqn:fullsequenceexample}
\phi : (D,\gamma,\tilde{\gamma}) = (D_1,\gamma_1,\tilde{\gamma}_1) \Move{h_1,\tilde{v}_1} (D_2,\gamma_2,\tilde{\gamma}_2) \Move{h_2,\tilde{v}_2} \cdots \Move{h_{n-1},\tilde{v}_{n-1}} (D_n,\gamma_n,\tilde{\gamma}_n) = (\ast,\ast,\tilde{v}).
\end{equation}
Such a sequence ending with $(\ast,\ast,\tilde{v})$ will be called a {\em full sequence of moves for $g$}.  We
will also discuss partial sequences of moves, which are defined in the obvious way.  Finally, we define
the {\em stabilizer product} of \eqref{eqn:fullsequenceexample} to be the product
$$(h_1)^{-1}_{\tilde{v}_1} \cdots (h_{n-1})^{-1}_{\tilde{v}_{n-1}} \cdot ((h_{n-1} \cdots h_1 \cdot g)^{-1})^{-1}_{\tilde{v}}  = (h_1)^{-1}_{\tilde{v}_1} \cdots (h_n)^{-1}_{\tilde{v}_n} \in \Gamma.$$
We will also discuss the stabilizer products of partial sequences of moves, which are again defined in the
obvious way (of course, in the stabilizer product of a partial sequence of moves, the $h_n$ term
is omitted).

\subsection{The dependence of the Armstrong construction on our choices}
\label{section:phiwelldefined}

Fix $g \in G$.  The goal of this section is to show that the element of $\Gamma$ obtained by applying the Armstrong
construction to $g$ is independent of all of our arbitrary choices.  Examining the construction, we see 
that the following is a complete list of such choices.

\begin{enumerate}
\item The path $\tilde{\gamma}$ from $\tilde{v}$ to $g(\tilde{v})$.
\item The nondegenerate disc $(D,\gamma)$ and the map $\phi:D \rightarrow X/G$ 
with $\phi_{\ast}(\gamma) = \pi_{\ast}(\tilde{\gamma})$.
\item The manner in which we collapsed $\gamma$ across $D$ to homotope it to the trivial loop.
\item The vertices $\tilde{y}$ involved in lifting two-dimensional collapses.
\item The group elements $h_i$ chosen during the liftings of one-dimensional collapses.
\end{enumerate}

\noindent
We will deal with each choice in turn.

\begin{remark}
The construction also depends on the choice of basepoint $\tilde{v} \in X^{(0)}$, but
that was fixed at the beginning so there is no need to show that the output
does not depend on it.  In fact, it is not hard to show that the output
of the Armstrong does not depend on $\tilde{v} \in X^{(0)}$ either.
\end{remark}

\begin{remark}[On the assumptions in the theorem]
The assumption that $X/G$ is simply connected was used during the Armstrong
construction.  The assumption that
$X$ is simply connected is used to show that the output of the Armstrong construction
does not depend on choice 1 and the assumption that $X/G$ is 2-connected is used to show
that the output of the Armstrong construction does not depend on choice 2.  These
assumptions are not used in the proof that the output of the Armstrong construction does
not depend on choices 3--5.
\end{remark}

\subsubsection{The output of the Armstrong construction does not depend on choices 3--5}

The proof will be by induction on the number $m$ of simplices (of any dimension) in the nondegenerate disc $D$.  In
fact, for use in our induction we will allow $(D,\gamma)$ to be a degenerate disc; as we observed in 
\S \ref{section:definephi}, the nondegeneracy of $(D,\gamma)$ was never used during the Armstrong construction.

The case $m=1$ is trivial; in that case, our expression for $\psi^{-1}(g)$ is simply $g_{\tilde{v}}$.  
Assume, therefore, that $m > 1$.  Consider two possible first moves
\begin{equation}
\label{eqn:firstmove1}
\phi : (D,\gamma,\tilde{\gamma}) \Move{h_1,\tilde{v}_1} (D_2,\gamma_2,\tilde{\gamma}_2)
\end{equation}
and
\begin{equation}
\label{eqn:firstmove2}
\phi : (D,\gamma,\tilde{\gamma}) \Move{h'_1,\tilde{v}_1'} (D_2',\gamma_2',\tilde{\gamma}_2').
\end{equation}
We will show that there exist full sequences of moves
\begin{equation}
\label{eqn:moves1}
\phi : (D,\gamma,\tilde{\gamma}) \Move{h_1,\tilde{v}_1} (D_2,\gamma_2,\tilde{\gamma}_2) \rightarrow \cdots 
\end{equation}
and
\begin{equation}
\label{eqn:moves2}
\phi : (D,\gamma,\tilde{\gamma}) \Move{h'_1,\tilde{v}_1'} (D_2',\gamma_2',\tilde{\gamma}_2') \rightarrow \cdots
\end{equation}
such that the stabilizer products of \eqref{eqn:moves1} and \eqref{eqn:moves2} are the same lift of $g$ to $\Gamma$.
This is enough to prove the claim.  Indeed, let 
\begin{equation}
\label{eqn:moves1prime}
\phi : (D,\gamma,\tilde{\gamma}) \Move{h_1,\tilde{v}_1} (D_2,\gamma_2,\tilde{\gamma}_2) \rightarrow \cdots 
\end{equation}
and 
\begin{equation}
\label{eqn:moves2prime}
\phi : (D,\gamma,\tilde{\gamma}) \Move{h'_1,\tilde{v}_1'} (D_2',\gamma_2',\tilde{\gamma}_2') \rightarrow \cdots
\end{equation}
be two arbitrary full sequences of moves across $D$ starting with \eqref{eqn:firstmove1} and \eqref{eqn:firstmove2}.
The subsequences of moves obtained by deleting the initial $(D,\gamma,\tilde{\gamma})$ from \eqref{eqn:moves1} and
\eqref{eqn:moves1prime} yield full sequences of moves across $D_2$.  The stabilizer products of these
subsequences give two lifts to $\Gamma$ of $h_1 g$, and
by our inductive hypothesis these two elements of $\Gamma$ are equal.  Thus the stabilizer products
of \eqref{eqn:moves1} and \eqref{eqn:moves1prime} are the same lift of $g$ to $\Gamma$.  
Similarly, the stabilizer products of \eqref{eqn:moves2} and \eqref{eqn:moves2prime}
are the same lift of $g$ to $\Gamma$.  
We conclude that the stabilizer products of \eqref{eqn:moves1prime} and \eqref{eqn:moves2prime} 
are the same lift of $g$ to $\Gamma$, as desired.

The heart of the proof will be the following three special cases.
\BeginCases
\begin{case}
\label{case:distinct}
The simplices of $D$ deleted in \eqref{eqn:firstmove1} and
\eqref{eqn:firstmove2} are distinct.  
\end{case}
\begin{case}
\label{case:line}
$D$ is a line
\end{case}
\begin{case}
\label{case:flower}
$D$ is a flower whose bloom consists of a single simplex.
\end{case}
These cases are enough to prove the claim.  Indeed, assume that we are not in Case \ref{case:distinct}, so 
the simplices of $D$ deleted in \eqref{eqn:firstmove1} and \eqref{eqn:firstmove2} are identical.  If 
there is some other first move
\begin{equation}
\label{eqn:firstmove3}
\phi : (D,\gamma,\tilde{\gamma}) \Move{h_1'',\tilde{v}_1''} (D_2'',\gamma_2'',\tilde{\gamma}_2'')
\end{equation}
which deletes a different simplex, then Case \ref{case:distinct} applies to \eqref{eqn:firstmove3} together with
either \eqref{eqn:firstmove1} or \eqref{eqn:firstmove2}.  Thus there exist full sequences of moves
\begin{align}
\phi : &(D,\gamma,\tilde{\gamma}) \Move{h_1,\tilde{v}_1} (D_2,\gamma_2,\tilde{\gamma}_2) \rightarrow \cdots \label{eqn:full1} \\
\phi : &(D,\gamma,\tilde{\gamma}) \Move{h_1'',\tilde{v}_1''} (D_2'',\gamma_2'',\tilde{\gamma}_2'') \rightarrow \cdots \label{eqn:full2} \\
\phi : &(D,\gamma,\tilde{\gamma}) \Move{h_1',\tilde{v}_1'} (D_2',\gamma_2',\tilde{\gamma}_2') \rightarrow \cdots \label{eqn:full3} \\
\phi : &(D,\gamma,\tilde{\gamma}) \Move{h_1'',\tilde{v}_1''} (D_2'',\gamma_2'',\tilde{\gamma}_2'') \rightarrow \cdots \label{eqn:full4}
\end{align}
such that the stabilizer products of \eqref{eqn:full1} and \eqref{eqn:full2} (resp.\ \eqref{eqn:full3} and \eqref{eqn:full4}) are
the same.  Using the inductive hypothesis like we did above, the stabilizer products of \eqref{eqn:full2} and \eqref{eqn:full4}
are the same.  We conclude that the stabilizer products of \eqref{eqn:full1} and \eqref{eqn:full3} are the same,
as desired.

If instead there is no move deleting a different simplex, then using Corollary 
\ref{corollary:lineflower} we can conclude that $D$ is a line (if \eqref{eqn:firstmove1}
and \eqref{eqn:firstmove2} are one-dimensional collapses) or a flower whose bloom consists of a single simplex 
(if \eqref{eqn:firstmove1} and \eqref{eqn:firstmove2} are two-dimensional collapses).  We can thus
apply either Case \ref{case:line} or Case \ref{case:flower}, and we are done.

The proofs of all three cases follow the same pattern.  Namely, for some $N \geq 2$ we construct partial sequences of moves
\begin{equation}
\label{eqn:firstmove1completed}
\phi : (D,\gamma,\tilde{\gamma}) \Move{h_1,\tilde{v}_1} (D_2,\gamma_2,\tilde{\gamma}_2) \Move{h_2,\tilde{v}_2} \cdots \Move{h_{N-1},\tilde{v}_{N-1}} (D_N,\gamma_N,\tilde{\gamma}_N)
\end{equation}
and
\begin{equation}
\label{eqn:firstmove2completed}
\phi : (D,\gamma,\tilde{\gamma}) \Move{h'_1,\tilde{v}_1'} (D_2',\gamma_2',\tilde{\gamma}_2')  \Move{h_2',\tilde{v}_2'} \cdots \Move{h_{N'-1}',\tilde{v}_{N'-1}'} (D_{N'}',\gamma_{N'}',\tilde{\gamma}_{N'}')
\end{equation}
with $(D_N,\gamma_N,\tilde{\gamma}_N) = (D_{N'}',\gamma_{N'}',\tilde{\gamma}_{N'}')$.  
We then verify that the stabilizer products
of \eqref{eqn:firstmove1completed} and \eqref{eqn:firstmove2completed} are equal; i.e.\ that 
$$(h_1)^{-1}_{\tilde{v}_1} \cdots (h_{N-1})^{-1}_{\tilde{v}_{N-1}} = (h_1')^{-1}_{\tilde{v}_1'} \cdots (h_{N'-1}')^{-1}_{\tilde{v}_{N'-1}}.$$
It will then follow that the desired full sequences of moves can be obtained by completing \eqref{eqn:firstmove1completed} and
\eqref{eqn:firstmove2completed} to full sequences of moves in the same way.

\Figure{figure:deletedisjoint}{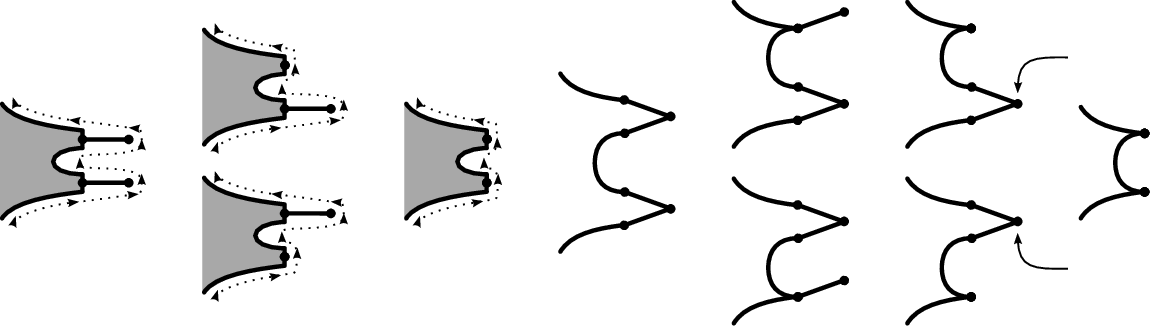}{a. Left is $(D_1,\gamma_1)$, middle top is $(D_2,\gamma_2)$, middle bottom is
$(D_2',\gamma_2')$, right is $(D_3,\gamma_3)$ \CaptionSpace b. The various lifts needed}

\Heading{Proof of Case 1 (distinct simplices deleted)}
There are three cases.
\begin{enumerate}
\item Both \eqref{eqn:firstmove1} and \eqref{eqn:firstmove2} correspond to two-dimensional collapses.
\item Both \eqref{eqn:firstmove1} and \eqref{eqn:firstmove2} correspond to one-dimensional collapses.
\item \eqref{eqn:firstmove1} and \eqref{eqn:firstmove2} correspond to different types of collapses.
\end{enumerate}
The argument is similar in all three cases; we will discuss the slightly more difficult case 2 and leave cases 1 and
3 to the reader.

Assume, therefore, that \eqref{eqn:firstmove1} and \eqref{eqn:firstmove2} correspond to one-dimensional collapses.  The
proof is illustrated in Figure \ref{figure:deletedisjoint}.  Let
$v_1$ and $v_1'$ be the vertices of $\gamma$ which map to $\pi(\tilde{v}_1)$ and $\pi(\tilde{v}_1')$.  Without loss
of generality, we can assume that $v_1'$ comes before $v_1$ in $\gamma$ (see Figure \ref{figure:deletedisjoint}.a).  Let
$(D_3,\gamma_3)$ be the result of performing a one-dimensional collapse to $(D_2,\gamma_2)$ at $v_1'$, or equivalently
the result of performing a one-dimensional collapse to $(D_2',\gamma_2)$ at $v_1$ (see Figure \ref{figure:deletedisjoint}.a).

Observe that swinging $\tilde{\gamma}$ around $\tilde{v}_1$ by $h_1$ does not affect $\tilde{v}_1'$, so 
we can lift $(\phi_3)_{\ast}(\gamma_3)$ to a path $\tilde{\gamma}_3$ in $X$ by
swinging $\tilde{\gamma}_2$ around $\tilde{v}_1'$ by $h_1'$ and then removing $\tilde{v}_1'$ from $\tilde{\gamma}_2$ 
(see Figure \ref{figure:deletedisjoint}.b).  Similarly,
swinging $\tilde{\gamma}$ around $\tilde{v}_1'$ by $h_1'$ moves $\tilde{v}_1$ to $h_1'(\tilde{v}_1)$, so
we can lift $(\phi_3)_{\ast}(\gamma_3)$ to $X$ by swinging 
$\tilde{\gamma}_2'$ around $h_1'(\tilde{v}_1)$ by $h_1' h_1 (h_1')^{-1}$ and
then removing $h_1'(\tilde{v}_1)$; the resulting path is equal to $\tilde{\gamma}_3$ (see 
Figure \ref{figure:deletedisjoint}.b).  

We thus have constructed sequences of moves
$$\phi : (D,\gamma,\tilde{\gamma}) \Move{h_1,\tilde{v}_1} (D_2,\gamma_2,\tilde{\gamma}_2) \Move{h_1',\tilde{v}_1'} (D_3,\gamma_3,\tilde{\gamma}_3)$$
and
$$\phi : (D,\gamma,\tilde{\gamma}) \Move{h_1',\tilde{v}_1'} (D_2',\gamma_2',\tilde{\gamma}_2') \Move{h_1' h_1 (h_1')^{-1},h_1'(\tilde{v}_1)} (D_3,\gamma_3,\tilde{\gamma}_3).$$
In $\Gamma$ we have
$$(h_1)_{\tilde{v}_1}^{-1} (h_1')_{\tilde{v}_1'}^{-1} = (h_1')_{\tilde{v}_1'}^{-1} (h_1')_{\tilde{v}_1'} (h_1)_{\tilde{v}_1}^{-1} (h_1')_{\tilde{v}_1'}^{-1} = (h_1')_{\tilde{v}_1'}^{-1} (h_1' h_1 (h_1')^{-1})_{h_1'(\tilde{v}_1)},$$
and the claim follows.  

\Heading{Proof of Case 2 ($\boldsymbol{D}$ a line)}
This proof is illustrated in Figure \ref{figure:deleteline}.  We will do the case that $D$ has at least $k$
vertices for $k \geq 3$.  The case that $D$ has 2 vertices (remember, it has more than 1) is similar, the key 
difference being that instead of the final ``swinging'' which occurs in our construction, use is made
of the final element (denoted $h_n$ in \S \ref{section:definephi}) in the stabilizer product of a full sequence of moves.

Let the last three vertices of our line be $w_{1}$, $w_{2}$, and $w_{3}$ (see Figure \ref{figure:deleteline}.a).  
As depicted in Figure \ref{figure:deleteline}.b, let the lift of the segment
$$\phi(w_{1}) - \phi(w_{2}) - \phi(w_3) - \phi(w_{2}) - \phi(w_{1})$$
of $\phi_{\ast}(\gamma)$ to $\tilde{\gamma}$ be
$$\tilde{w}_{1} - \tilde{w}_{2} - \tilde{w}_{3} - \tilde{w}_{2}' - \tilde{w}_{1}';$$
Observe that in \eqref{eqn:firstmove1} and \eqref{eqn:firstmove2}, we must have $v_1 = v_1' = \tilde{w}_3$.  We therefore
have $(D_2,\gamma_2) = (D_2',\gamma_2')$ (see Figure \ref{figure:deleteline}.a).  Let $(D_3,\gamma_3)$ be
as in Figure \ref{figure:deleteline}.a.

\Figure{figure:deleteline}{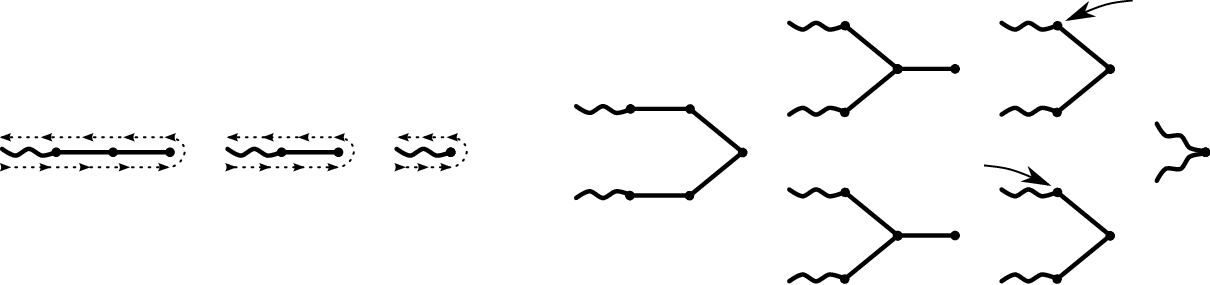}{a. Left is $(D,\gamma)$, middle is $(D_2,\gamma_2) = (D_2',\gamma_2')$, right
is $(D_3,\gamma_3)$ \CaptionSpace b. Lifts to $X$}

Next, observe that
$$h_1(\tilde{w}_{2}') = h_1'(\tilde{w}_{2}') = \tilde{w}_{2}.$$
This implies that $h_1' = \mu h_1$ for some $\mu \in G_{\{\tilde{w}_3,\tilde{w}_{2}\}}$.  Observe
that $\tilde{\gamma}_2$ and $\tilde{\gamma}_2'$ are as depicted in Figure \ref{figure:deleteline}.b.
Let $h_2 \in G_{\tilde{w}_{2}}$ be such that
$$h_2(h_1(\tilde{w}_{1}')) = \tilde{w}_{1}.$$
We can therefore obtain a lift $\tilde{\gamma}_3$ of $\phi_{\ast}(\gamma_3)$ to $X$ by swinging 
$\tilde{\gamma}_2$ around the pivot point $\tilde{w}_{2}$ by $h_2$ and then deleting $\tilde{w}_{2}$
from it.  We thus have a sequence of moves
\begin{equation}
\label{eqn:partialline1}
\phi : (D,\gamma,\tilde{\gamma}) \Move{h_1,\tilde{w}_3} (D_2,\gamma_2,\tilde{\gamma}_2) \Move{h_2,\tilde{w}_{2}} (D_3,\gamma_3,\tilde{\gamma}_3)
\end{equation}
as depicted in Figure \ref{figure:deleteline}.b.

Since $\mu,h_2 \in G_{\tilde{w}_{2}}$, it follows that
$h_2 \mu^{-1} \in G_{\tilde{w}_{2}}$.  Also, we have 
$$(h_2 \mu^{-1})(h_1'(\tilde{w}_{1}')) = h_2 \mu^{-1} \mu h_1 (\tilde{w}_{1}') = h_2(h_1(\tilde{w}_{1}')) = \tilde{w}_{1}.$$
We can thus obtain a lift $\tilde{\gamma}_3'$ of $\phi_{\ast}(\gamma_3')$
to $X$ by swinging $\tilde{\gamma}_2'$ around the pivot point $\tilde{w}_{2}$ by $h_2 \mu^{-1}$
and then deleting $\tilde{w}_{2}$ from it.  Moreover, it is clear that $\tilde{\gamma}_3' = \tilde{\gamma}_3$.  We thus
have a sequence of moves
\begin{equation}
\label{eqn:partialline2}
\phi : (D,\gamma,\tilde{\gamma}) \Move{h_1',\tilde{w}_3} (D_2',\gamma_2',\tilde{\gamma}_2') \Move{h_2 \mu^{-1},\tilde{w}_{2}} (D_3,\gamma_3,\tilde{\gamma}_3).
\end{equation}
Since in $\Gamma$ we have
$$(h_1')^{-1}_{\tilde{w}_3} (h_2 \mu^{-1})^{-1}_{\tilde{w}_{2}} = (h_1)^{-1}_{\tilde{w}_3} (\mu)^{-1}_{\tilde{w}_3} (\mu)_{\tilde{w}_{2}} (h_2)^{-1}_{\tilde{w}_{2}} = (h_1)^{-1}_{\tilde{w}_3} (\mu)^{-1}_{\tilde{w}_3} (\mu)_{\tilde{w}_{3}} (h_2)^{-1}_{\tilde{w}_{2}} = (h_1)^{-1}_{\tilde{w}_3} (h_2)^{-1}_{\tilde{w}_{2}},$$
the claim follows.

\Heading{Proof of Case 3 ($\boldsymbol{D}$ a flower)}
We will assume that the stem of our flower has at least two vertices (i.e.\ that there is more
to the flower than the bloom -- the bloom and the stem share one vertex); the case where the stem
consists of a single vertex is similar.  As in the top portion 
of Figure \ref{figure:deleteflower}.a, let $w_1 - w_2$ be the last edge of the stem of $D$ and let
the path around the bloom of $D$ be $w_2 - x - y - w_2$.  There are two cases.

\BeginSubcases
\begin{subcase}
The edges deleted by \eqref{eqn:firstmove1} and \eqref{eqn:firstmove2} are equal.
\end{subcase}

\Figure{figure:deleteflower}{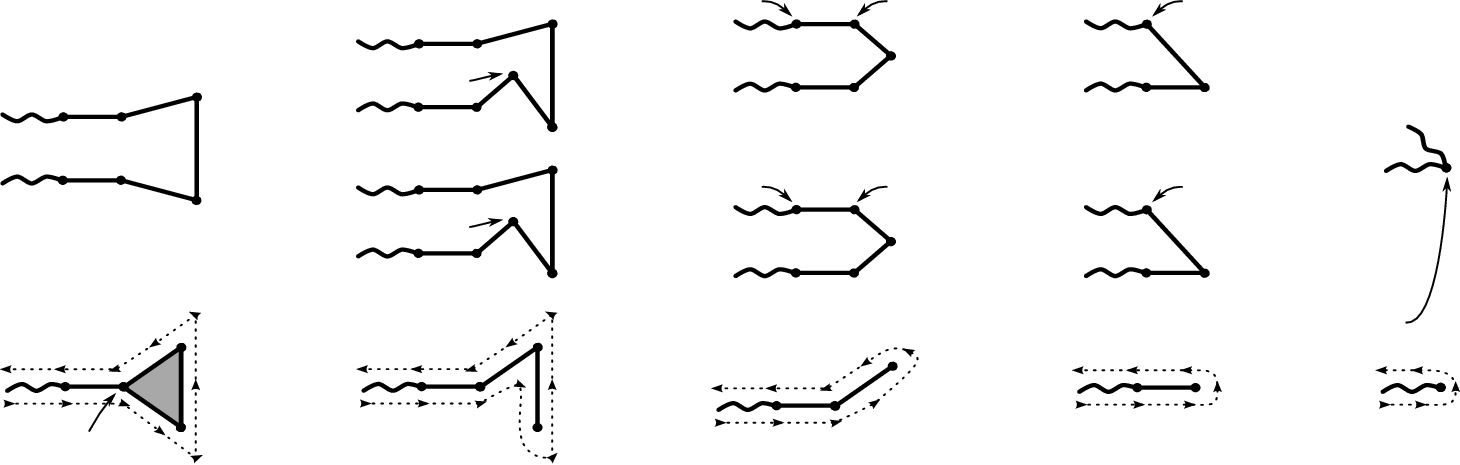}{a. Bottom is $(D,\gamma)$, top is $\tilde{\gamma}$
\CaptionSpace b--e. Bottoms are $(D_2,\gamma_2) = (D_2',\gamma_2')$ through $(D_5,\gamma_5) = (D_5',\gamma_5')$, tops
are lifts to $X$}

This is similar to the case of the line, so our exposition will be slightly terse.  We will do the case that
\eqref{eqn:firstmove1} and \eqref{eqn:firstmove2} delete the edge $w_2 - x$; the other two cases
are similar.  Observe that $(D_2,\gamma_2) = (D_2',\gamma_2')$ and
that $(D_2,\gamma_2)$ is as depicted on the bottom of Figure \ref{figure:deleteflower}.b.  Also, since \eqref{eqn:firstmove1}
and \eqref{eqn:firstmove2} are both two-dimensional collapses, it follows that $h_1 = h_1' = 1$.  Let $(D_i,\gamma_i)$ for 
$3 \leq i \leq 5$ be as in the bottom portions of Figures \ref{figure:deleteflower}.c--e.

Let the portion of $\tilde{\gamma}_1$ corresponding to the last two vertices of the stem plus the 
path around the bloom be
$$\tilde{w}_1 - \tilde{w}_2 - \tilde{x} - \tilde{y} - \tilde{w}_2' - \tilde{w}_1'$$
and let the corresponding portion of $\tilde{\gamma}_2$ be
\begin{equation}
\label{eqn:flowergamma2lift}
\tilde{w}_1 - \tilde{w}_2 - \tilde{y}' - \tilde{x} - \tilde{y} - \tilde{w}_2' - \tilde{w}_1';
\end{equation}
see the tops of Figures \ref{figure:deleteflower}.a,b.  Then as shown in Figures \ref{figure:deleteflower}.a--e, 
there exists $h_2 \in G_{\tilde{x}}$, $h_3 \in G_{\tilde{y}'}$, and
$h_4 \in G_{\tilde{w}_2}$ such that we have a sequence of moves
$$\phi : (D,\gamma,\tilde{\gamma}) \Move{h_1,\tilde{v}_1} (D_2,\gamma_2,\tilde{\gamma}_2)
\Move{h_2,\tilde{x}} (D_3,\gamma_3,\tilde{\gamma}_3) \Move{h_3,\tilde{y}'} (D_4,\gamma_4,\tilde{\gamma}_4)
\Move{h_4,\tilde{w}_2} (D_5,\gamma_5,\tilde{\gamma}_5).$$

Moreover, as in the case of the line, there must exist some $\mu \in G_{\{\tilde{w}_2,\tilde{x}\}}$ such that
the lift of the portion of $\tilde{\gamma}_2'$ corresponding to \eqref{eqn:flowergamma2lift} is
$$\tilde{w}_1 - \tilde{w}_2 - \mu(\tilde{y}') - \tilde{x} - \tilde{y} - \tilde{w}_2' - \tilde{w}_1'.$$
Sine $\mu \in G_{\{\tilde{w}_2,\tilde{x}\}}$, we have 
$\mu h_2 \in G_{\tilde{x}}$ and $h_4 \mu^{-1} \in G_{\tilde{w}_2}$.  Thus
as shown in Figures \ref{figure:deleteflower}.a--e, there exists a sequence of moves
$$\phi : (D,\gamma,\tilde{\gamma}) \Move{h_1',\tilde{v}_1'} (D_2',\gamma_2',\tilde{\gamma}_2')
\Move{\mu h_2,\tilde{x}} (D_3',\gamma_3',\tilde{\gamma}_3') \Move{\mu h_3 \mu^{-1},\mu(\tilde{y}')} 
(D_4',\gamma_4',\tilde{\gamma}_4') \Move{h_4 \mu^{-1},\tilde{w}_2} (D_5',\gamma_5',\tilde{\gamma}_5').$$
Also, we have $(D_5',\gamma_5',\tilde{\gamma}_5') = (D_5,\gamma_5,\tilde{\gamma}_5)$.  As we observed above,
we have $h_1 = h_1' = 1$, so the proof is completed by observing that in $\Gamma$ we have
\begin{align*}
(\mu h_2)^{-1}_{\tilde{x}} (\mu h_3 \mu^{-1})^{-1}_{\mu(\tilde{y}')} (h_4 \mu^{-1})^{-1}_{\tilde{w}_2}
&= (h_2)^{-1}_{\tilde{x}} (\mu)^{-1}_{\tilde{x}} (\mu)_{\tilde{x}} (h_3)^{-1}_{\tilde{y}'} (\mu)^{-1}_{\tilde{x}} (\mu)_{\tilde{w}_2} (h_4)^{-1}_{\tilde{w}_2}\\
&= (h_2)^{-1}_{\tilde{x}} (h_3)^{-1}_{\tilde{y}'} (\mu)^{-1}_{\tilde{x}} (\mu)_{\tilde{x}} (h_4)^{-1}_{\tilde{w}_2}\\
&= (h_2)^{-1}_{\tilde{x}} (h_3)^{-1}_{\tilde{y}'} (h_4)^{-1}_{\tilde{w}_2}
\end{align*}

\begin{subcase}
The edges deleted by \eqref{eqn:firstmove1} and \eqref{eqn:firstmove2} are different.
\end{subcase}

\Figure{figure:deleteflower2}{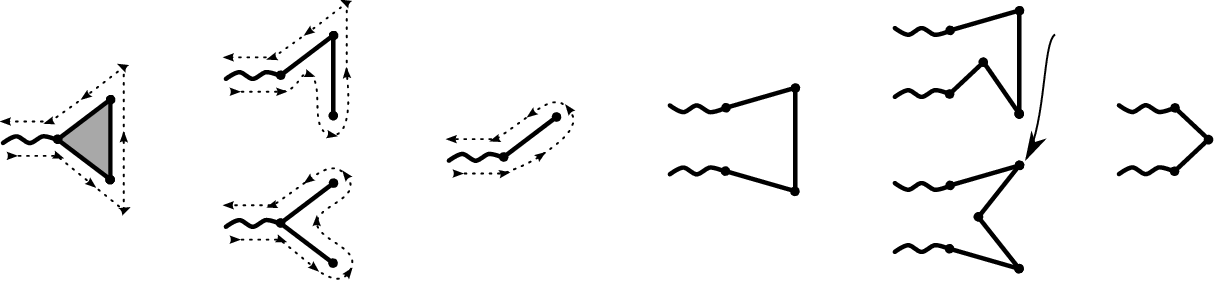}{a. Left is $\gamma$, middle top is $\gamma_2$, middle bottom is $\gamma_2'$, 
right is $\gamma_3$
\CaptionSpace
b. Various lifts needed.}

We will discuss the case that \eqref{eqn:firstmove1} deletes the edge $w_2 - x$ and \eqref{eqn:firstmove2}
deletes the edge $x - y$; the other cases are similar.  Again, the details are similar to what we have
already done, so we only sketch the argument.

Observe that $(D_2,\gamma_2)$ and $(D_2',\gamma_2')$ 
are as depicted in the top and bottom of the middle of Figure \ref{figure:deleteflower2}.a, respectively.  
Also, since \eqref{eqn:firstmove1}
and \eqref{eqn:firstmove2} are both two-dimensional collapses, it follows that $h_1 = h_1' = 1$.  Let $(D_3,\gamma_3)$ 
be as in the right hand portion of Figure \ref{figure:deleteflower2}.a.

As in Figure \ref{figure:deleteflower2}.b, let the portions of $\tilde{\gamma}_2$ and $\tilde{\gamma}_2'$ corresponding to the 
last vertex of the stem plus the path around the bloom be
$$\tilde{w}_2 - \tilde{y}' - \tilde{x} - \tilde{y} - \tilde{w}_2'$$
and
$$\tilde{w}_2 - \tilde{x} - \tilde{w}_2'' - \tilde{y} - \tilde{w}_2',$$
respectively (we will discuss the identities $\tilde{w}_2'' = h^{-1}(\tilde{w}_2)$ and
$\tilde{y} = h^{-1}(\tilde{y}')$ in this figure shortly).  Also, let
$$\phi|_{D_2} : (D_2,\gamma_2,\tilde{\gamma}_2) \Move{h_2,\tilde{x}} (D_3,\gamma_3,\tilde{\gamma}_3)$$
be an arbitrary lift of a one-dimensional collapse.  Thus $h_2 \in G_{\tilde{x}}$ and $h_2(\tilde{y}) = \tilde{y}'$.
Using the previous case, we can modify $\tilde{w}_2''$ to {\em any} legal vertex.  In particular,
we can assume that $\tilde{w}_2'' = h_2^{-1}(\tilde{w}_2)$ (see Figure \ref{figure:deleteflower2}.b).  Thus
$$\phi|_{D_2'} : (D_2',\gamma_2',\tilde{\gamma}_2') \Move{h_2,\tilde{x}} (D_3,\gamma_3,\tilde{\gamma}_2)$$
is a lift of a one-dimensional collapse.  Summing up, we have sequences of moves 
$$\phi : (D,\gamma,\tilde{\gamma}) = (D_1,\gamma_1,\tilde{\gamma}_1) \Move{h_1,\tilde{v}_1} (D_2,\gamma_2,\tilde{\gamma}_2)
\Move{h_2,\tilde{x}} (D_3,\gamma_3,\tilde{\gamma}_3)$$
and
$$\phi : (D,\gamma,\tilde{\gamma}) = (D_1',\gamma_1',\tilde{\gamma}_1') \Move{h_1',\tilde{v}_1'} (D_2',\gamma_2',\tilde{\gamma}_2') 
\Move{h_2,\tilde{x}} (D_3,\gamma_3,\tilde{\gamma}_3).$$ 
Since
$$(h_1)_{\tilde{v}_1}^{-1} (h_2)_{\tilde{x}}^{-1} = 1 \cdot (h_2)_{\tilde{x}}^{-1} = (h_1')_{\tilde{v}_1}^{-1} (h_2)_{\tilde{x}}^{-1},$$
the result follows.

\subsubsection{The output of the Armstrong construction does not depend on choice 2 (the degenerate disc $\boldsymbol{(D,\gamma)}$).}

Recall that we have fixed some $g \in G$.  By the previous section, the Armstrong construction associates a well-defined
element of $\Gamma$ to the following data.
\begin{itemize}
\item A simplicial path $\tilde{\gamma}$ in $X$ from $\tilde{v}$ to $g(\tilde{v})$.
\item A nondegenerate disc $(D,\gamma)$.
\item A simplicial map $\phi : D \rightarrow X/G$ such that $\phi_{\ast}(\gamma) = \pi_{\ast}(\tilde{\gamma})$.
\end{itemize}
In this section, we show that output of the Armstrong construction is independent of $(D,\gamma)$ and
$\phi$.

Fix a path $\tilde{\gamma}$, a nondegenerate disc $(D,\gamma)$, and a map $\phi : D \rightarrow X/G$ as above.
Since $X/G$ is $2$-connected, Theorem \ref{theorem:elementaryequiv} and the fact that any two triangulations
of a closed $2$-disc have a common subdivision imply that it is enough to show
that the output of the Armstrong construction is invariant under elementary subdivisions and elementary
pushes.  The arguments for these are similar; we will give the details for elementary pushes and
leave the case of elementary subdivisions to the reader.

Let $y$ be a vertex of $D$ that does not lie on $\gamma$.  Set $z = \phi(y)$ and let
$z' \in X/G$ be any vertex such that for all simplices $s$ of $\Star_{D}(y)$, the set
$\phi(s) \cup z'$ is a simplex of $X/G$.  Define $\phi' : D \rightarrow X/G$ to equal $\phi$ on all 
vertices except for $y$, where it equals $z'$.  
Our goal is to show that there exist full sequences of moves
\begin{equation}
\label{eqn:deleteflowerfull1}
\phi : (D,\gamma,\tilde{\gamma}) \rightarrow \cdots
\end{equation}
and
\begin{equation}
\label{eqn:deleteflowerfull2}
\phi' : (D,\gamma,\tilde{\gamma}) \rightarrow \cdots
\end{equation}
such that the stabilizer products of \eqref{eqn:deleteflowerfull1} and \eqref{eqn:deleteflowerfull2} are equal.

\Figure{figure:bigflowercollapse}{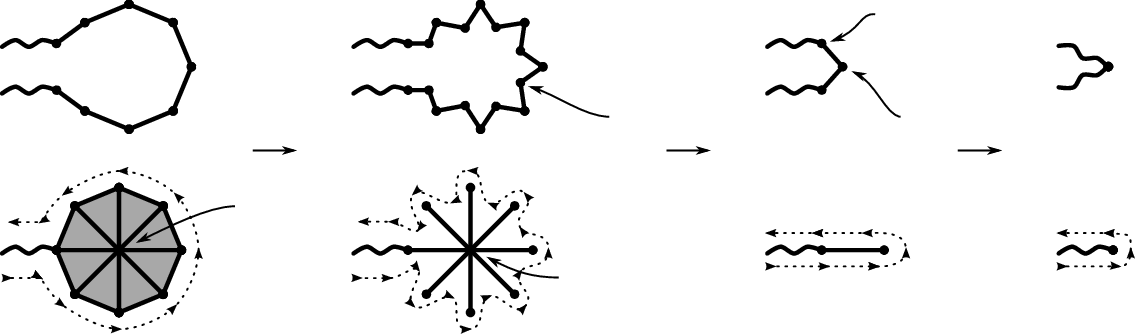}{a. Bottom is $(D,\gamma)$, top is $\tilde{\gamma}$
\CaptionSpace b. Bottom is $(E,\delta)$, top is $\tilde{\delta}$ and $\tilde{\delta}'$
\CaptionSpace c. Bottom is $(F,\epsilon)$, top is $\tilde{\epsilon}$ and $\tilde{\epsilon}'$
\CaptionSpace d. Bottom is $(C,\eta)$, top is $\tilde{\eta}$}

By Corollary \ref{corollary:lineflower}, we can assume that $(D,\gamma)$ is a flower whose
bloom is $\Star_{D}(y)$ (see Figure \ref{figure:bigflowercollapse}.a).  Let the images
under $\phi$ and $\phi'$ of the portion of $\gamma$ that goes around the bloom (these are
equal) be
$$x_0 - x_1 - \cdots - x_k = x_0,$$
and let the corresponding subpath of $\tilde{\gamma}$ be
\begin{equation}
\label{eqn:flowercollapsebloom}
\tilde{x}_0 - \tilde{x}_1 - \cdots - \tilde{x}_k.
\end{equation}
For $0 \leq i < k$, let $\tilde{z}_i$ be a lift of $z$ to $X$ and $\tilde{z}_i'$ be a lift
of $z'$ to $X$ such that $\{\tilde{x}_i, \tilde{x}_{i+1}, \tilde{z}_i, \tilde{z}_i'\}$ is a simplex
of $X$.  Define $(E,\delta)$ to be as depicted in Figure \ref{figure:bigflowercollapse}.b.  Let
$\tilde{\delta}$ and $\tilde{\delta}'$ to be $\tilde{\gamma}$ with the subpath 
\eqref{eqn:flowercollapsebloom}
replaced with the paths 
$$\tilde{x}_0 - \tilde{z}_0 - \tilde{x}_1 - \tilde{z}_1 - \tilde{x}_2 - \cdots - \tilde{z}_{k-1} - \tilde{x}_k$$
and
$$\tilde{x}_0 - \tilde{z}_0' - \tilde{x}_1 - \tilde{z}_1' - \tilde{x}_2 - \cdots - \tilde{z}_{k-1}' - \tilde{x}_k,$$
respectively (see Figure \ref{figure:bigflowercollapse}.b).  
Observe that there is a sequence of lifted $2$-dimensional collapses which converts
$(D,\gamma,\tilde{\gamma})$ to $(E,\delta,\tilde{\delta})$, and similarly there is a sequence of lifted
$2$-dimensional collapses which converts $(D,\gamma,\tilde{\gamma})$ to $(E,\delta,\tilde{\delta}')$.

Now let $(F,\epsilon)$ equal the degenerate disc depicted in Figure \ref{figure:bigflowercollapse}.c; 
it is obtained by performing one-dimensional collapses to the ``spokes'' of $E$.  
We lift these one-dimensional collapses as follows.  First, we choose $h_1 \in G_{\tilde{x}_1}$
such that $h_1(\tilde{z}_1) = \tilde{z}_0$ and $h_1(\tilde{z}_1') = \tilde{z}_0'$; such a choice is possible
since $\{\tilde{x}_1,\tilde{z}_0,\tilde{z}_0'\}$ and $\{\tilde{x}_1,\tilde{z}_1,\tilde{z}_1'\}$ are both
simplices of $X$.  We then swing $\tilde{\delta}$ and $\tilde{\delta}'$ around the pivot point $\tilde{x}_1$ by $h_1$ and
delete $\tilde{x}_1$.  Next, we choose $h_2 \in G_{h_1(\tilde{x}_1)}$ such that
$h_2(h_1(\tilde{z}_2)) = \tilde{z}_0$ and $h_2(h_1(\tilde{z}_2')) = \tilde{z}_0'$; such a choice
is possible since both
$$h_1(\{\tilde{x}_2,\tilde{z}_2,\tilde{z}_2'\}) = \{h_1(\tilde{x}_2),h_1(\tilde{z}_2),h_1(\tilde{z}_2')\} \quad \text{and}
\quad h_1(\{\tilde{x}_2,\tilde{z}_1,\tilde{z}_1'\}) = \{h_1(\tilde{x}_2),\tilde{z}_0,\tilde{z}_0'\}$$
are simplices of $X/G$.  We then swing $\tilde{\delta}$ and $\tilde{\delta}'$ around the pivot point $h_1(\tilde{x}_2)$ by
$h_2$ and delete $h_1(\tilde{x}_2)$.  Continuing in this manner, we obtain lifts $\tilde{\epsilon}$
and $\tilde{\epsilon}'$ of $\epsilon$ such that there are sequences
of one-dimensional collapses
\begin{equation}
\label{eqn:collapsespokes1}
\phi|_{E} : (E,\delta,\tilde{\delta}) \rightarrow \cdots \rightarrow (F,\epsilon,\tilde{\epsilon})
\end{equation}
and
\begin{equation}
\label{eqn:collapsespokes2}
\phi'|_{E} : (E,\delta,\tilde{\delta}') \rightarrow \cdots \rightarrow (F,\epsilon,\tilde{\epsilon}')
\end{equation}
satisfying the following two conditions.
\begin{itemize}
\item In \eqref{eqn:collapsespokes1} and \eqref{eqn:collapsespokes2}, the paths are swung around
the {\em same} pivot points by the {\em same} group elements.
\item The paths $\tilde{\epsilon}$ and $\tilde{\epsilon}'$ are as depicted in 
Figure \ref{figure:bigflowercollapse}.c.
\end{itemize} 

Observe that the sets
$$\{h_{k-1} \cdots h_1(\tilde{x}_{k}),h_{k-1} \cdots h_1(\tilde{z}_{k-1}), h_{k-1} \cdots h_1(\tilde{z}_{k-1}')\} \quad \text{and} \quad
\{\tilde{x}_0,h_{k-1} \cdots h_1(\tilde{z}_{k-1}), h_{k-1} \cdots h_1(\tilde{z}_{k-1}')\}$$
are simplices of $X/G$.  Thus there exists some 
$$h_k \in G_{\{h_{k-1} \cdots h_1(\tilde{z}_{k-1}), h_{k-1} \cdots h_1(\tilde{z}_{k-1}')\}}$$
such that $h_k h_{k-1} \cdots h_1(\tilde{x}_{k}) = \tilde{x}_0$.  Let $(C,\eta)$ be the degenerate
disc depicted in Figure \ref{figure:bigflowercollapse}.d and let $\tilde{\eta}$ be the
result of swinging $\tilde{\epsilon}$ around $h_{k-1} \cdots h_1(\tilde{z}_{k-1})$ by $h_k$ and then
deleting $h_{k-1} \cdots h_1(\tilde{z}_{k-1})$ (see Figure \ref{figure:bigflowercollapse}.d).  Observe that
$\tilde{\eta}$ also equals the result of swinging $\tilde{\epsilon}'$ around
$h_{k-1} \cdots h_1(\tilde{z}_{k-1}')$ by $h_k$ and then deleting $h_{k-1} \cdots h_1(\tilde{z}_{k-1}')$.

Summing up, we have constructed sequences of moves
\begin{equation}
\label{eqn:finalequation1}
\phi : (D,\gamma,\tilde{\gamma}) \rightarrow \cdots \rightarrow (E,\delta,\tilde{\delta}) \rightarrow \cdots \rightarrow (F,\epsilon,\tilde{\epsilon}) \Move{h_k,h_{k-1} \cdots h_1(\tilde{z}_{k-1})} (C,\eta,\tilde{\eta})
\end{equation}
and
\begin{equation}
\label{eqn:finalequation2}
\phi' : (D,\gamma,\tilde{\gamma}) \rightarrow \cdots \rightarrow (E,\delta,\tilde{\delta}') \rightarrow \cdots \rightarrow (F,\epsilon,\tilde{\epsilon}') \Move{h_k,h_{k-1} \cdots h_1(\tilde{z}_{k-1}')} (C,\eta,\tilde{\eta})
\end{equation}
Moreover, in both of these sequences the pivot points and group elements are all identical except for the final ones.  Since
$\tilde{z}_{k-1}$ and $\tilde{z}_{k-1}'$ are joined by an edge, we have an edge relation
$$(h_k)_{h_{k-1} \cdots h_1(\tilde{z}_{k-1})} = (h_k)_{h_{k-1} \cdots h_1(\tilde{z}_{k-1}')}$$
in $\Gamma$.  Thus the stabilizer products of \eqref{eqn:finalequation1} and
\eqref{eqn:finalequation2} are equal, and we are done.

\subsubsection{The output of the Armstrong construction does not depend on choice 1 (the path $\boldsymbol{\tilde{\gamma}}$ from $\boldsymbol{\tilde{v}}$ to $\boldsymbol{g(\tilde{v})}$)}

By the previous two sections, the output of the Armstrong construction only depends on the path
$\tilde{\gamma}$ from $\tilde{v}$ to $g(\tilde{v})$.  We wish to show that in fact it is independent
of this path.

\Figure{figure:lineindependent}{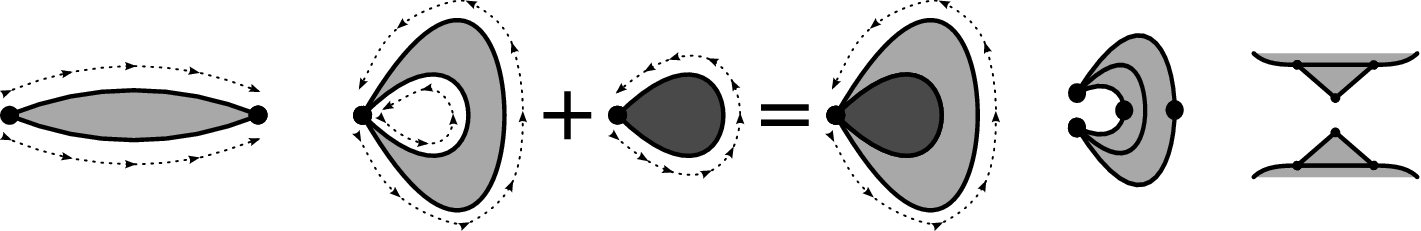}{a. $(Q,\ast)$
\CaptionSpace b. Gluing $A$ and $D$ to get $D'$
\CaptionSpace c. Gluing $\ast$ to $x$ does not result in a simplicial complex
\CaptionSpace d. Gluing top to bottom does not result in a simplicial complex}

Consider any two simplicial paths $\tilde{\gamma}$ and $\tilde{\gamma}'$ from $\tilde{v}$ to $g(\tilde{v})$.
We will construct nondegenerate discs $(D,\gamma)$ and $(D',\gamma')$ together with simplicial maps
$\phi : D \rightarrow X/G$ and $\phi' : D' \rightarrow X/G$ such that $\phi_{\ast}(\gamma) = \pi_{\ast}(\tilde{\gamma})$
and $\phi'_{\ast}(\gamma') = \pi_{\ast}(\tilde{\gamma}')$.  Additionally, we will construct full sequences
of moves
\begin{equation}
\label{eqn:pathindepfull1}
\phi : (D,\gamma,\tilde{\gamma}) \rightarrow \cdots
\end{equation}
and
\begin{equation}
\label{eqn:pathindepfull2}
\phi : (D',\gamma',\tilde{\gamma}') \rightarrow \cdots
\end{equation}
such that the stabilizer products of \eqref{eqn:pathindepfull1} and \eqref{eqn:pathindepfull2} are equal.  This
will imply the desired result.

Since $X$ is simply connected, $\tilde{\gamma}$ and $\tilde{\gamma}'$ are homotopic fixing the endpoints.  
This implies that there exists a based simplicial complex $(Q,\ast)$ such that $Q$ is homeomorphic to a closed
$2$-disc and $\ast \in \partial Q$ together with a map $\tilde{\rho} : Q \rightarrow X$ with the following property.  
There is a vertex $x$ on $\partial Q$ such that
if $\epsilon$ and $\epsilon'$ are the two embedded paths in $\partial Q$ from $\ast$ to $x$, then
$\tilde{\rho}_{\ast}(\epsilon) = \tilde{\gamma}$ and $\tilde{\rho}_{\ast}(\epsilon') = \tilde{\gamma}'$
(see Figure \ref{figure:lineindependent}.a).  Subdividing $Q$ if necessary away from $\partial Q$, we can
glue $x$ and $\ast$ together to get a simplicial complex $A$ (see Figure \ref{figure:lineindependent}.c for
an example of why we may need to subdivide $Q$ first).  Observing that $\pi(\tilde{\rho}(x)) = \pi(\tilde{\rho}(\ast)) = v$,
let $\rho : A \rightarrow X/G$ be the map induced by
by $\pi \circ \tilde{\rho} : Q \rightarrow X/G$ (see Figure \ref{figure:lineindependent}.b) and let the
projections of the paths $\epsilon$ and $\epsilon'$ in $Q$ to $A$ be $\overline{\epsilon}$ and $\overline{\epsilon}'$,
respectively.

Choose a nondegenerate disc $(D,\gamma)$ and a map $\phi : D \rightarrow X/G$ such that 
$\phi_{\ast}(\gamma) = \pi_{\ast}(\tilde{\gamma})$.
Subdividing $D$ away from $\gamma$ and $A$ away from $\overline{\epsilon}$ and $\overline{\epsilon}'$ if necessary,
glue $\overline{\epsilon} \subset A$ to $\gamma \subset D$ to obtain a simplicial complex 
$D'$ which is homeomorphic to a closed two-disc (see Figure \ref{figure:lineindependent}.b; also see
Figure \ref{figure:lineindependent}.d for an example of why we may need to subdivide).
Let $\gamma'$ be the loop around the boundary of $D'$ corresponding to $\overline{\epsilon}'$.
Thus $(D',\gamma')$ is another nondegenerate disc.  The maps $\phi$ and $\rho$ induce a map 
$\phi' : D' \rightarrow X/G$ such that
$\phi'_{\ast}(\gamma') = \pi_{\ast}(\tilde{\gamma}')$.  Now, by Corollary \ref{corollary:lineflower},
there is a sequence of collapses which converts $(D',\gamma')$ into $(D,\gamma)$.  Moreover,
by construction there is a partial sequence of moves
\begin{equation}
\label{eqn:pathindeppartial}
\phi' : (D',\gamma',\tilde{\gamma}') \rightarrow \cdots \rightarrow (D,\gamma,\tilde{\gamma})
\end{equation}
in which no swinging occurs (and, in particular, whose stabilizer product is trivial).  Letting
\begin{equation}
\label{eqn:pathindepfull1blah}
\phi : (D,\gamma,\tilde{\gamma}) \rightarrow \cdots
\end{equation}
be any full sequence of moves across $(D,\gamma,\tilde{\gamma})$, we can append \eqref{eqn:pathindeppartial}
to the beginning of \eqref{eqn:pathindepfull1blah} to obtain a full sequence of
moves
\begin{equation}
\label{eqn:pathindepfull2blah}
\phi' : (D',\gamma',\tilde{\gamma}') \rightarrow \cdots \rightarrow (D,\gamma,\tilde{\gamma}) \rightarrow \cdots.
\end{equation}
By construction, the stabilizer products of \eqref{eqn:pathindepfull1blah} and \eqref{eqn:pathindepfull2blah} are equal, and
we are done.

\subsection{Proof that the Armstrong construction defines a surjective homomorphism}
\label{section:phihomo}

In \S \ref{section:phiwelldefined}, we showed that the Armstrong construction gives a well-defined
map $\psi^{-1} : G \rightarrow \Gamma$.  It is clear that $\psi^{-1}$ is a homomorphism.  We must 
check that it is surjective and that 
$\psi \circ \psi^{-1} = 1$.  Observe that $\Gamma$ is generated by the set
$$S = \{\text{$g_{\tilde{x}}$ $|$ $\tilde{x} \in X^{(0)}$ and $g \in G_{\tilde{x}}$}\}.$$
Consider $g_{\tilde{x}} \in S$.  It is enough to show that $\psi^{-1}(g) = g_{\tilde{x}}$.

Let $\tilde{\gamma}'$ be a path in $X$ from 
$\tilde{v}$ to $\tilde{x}$ and let $(\tilde{\gamma}')^{-1}$ be $\tilde{\gamma}'$ traversed in reverse order.  
We then obtain a path $\tilde{\gamma}$ in $X$ from $\tilde{v}$ to $g(\tilde{v})$ by concatenating 
$\tilde{\gamma}'$ with $g((\tilde{\gamma}')^{-1})$.  Observe that $\pi_{\ast}(\tilde{\gamma})$ is
the image of the boundary of a line.  We can construct a full sequence of moves across this 
line by first swinging $\tilde{\gamma}$ around $\tilde{x}$ by $g^{-1}$ and then
doing a sequence of one-dimensional collapses without any additional swinging.  This implies that 
$\psi^{-1}(g) = g_{\tilde{x}}$, and we are done.

\noindent
Department of Mathematics\\
Rice University, MS 136\\
6100 Main St.\\
Houston, TX 77005\\
E-mail: {\tt andyp@math.rice.edu}
\medskip

\end{document}